\documentclass[12pt]{article}
\usepackage{amsmath}
\usepackage{amsfonts}
\usepackage{amssymb}

\title{Mixed identities, hereditarily separated actions and oscillation }

\author{A. Ivanov and R. Zarzycki}

\newtheorem{theorem}{Theorem}[section]
\newtheorem{proposition}[theorem]{Proposition}
\newtheorem{corollary}[theorem]{Corollary}
\newtheorem{lemma}[theorem]{Lemma}
\newtheorem{definition}[theorem]{Definition}
\newtheorem{example}[theorem]{Example}
\newtheorem{remark}[theorem]{Remark}

\begin{document} 
\topmargin = 12pt
\textheight = 630pt 
\footskip = 39pt 

\maketitle

\begin{quote}
{\bf Abstract} 
Given a topological $G$-space we consider equations with constants over $G$. 
In particular we formulate some very general conditions on words with constants $w(\bar{y},\bar{g})$ over $G$ which guarantee that the inequality $w(\bar{y},\bar{g})\neq 1$ has a solution in $G$. 
These results are illustrated in some typical situations. 
In particular standard actions of Thompson's group $F$ and branch groups are considered. 
\end{quote}

\section{Introduction} 

Let $\mathbb{F}_n$ denote the free group of rank $n$. 
For $w\in \mathbb{F}_n \ast G$ the equality $w=1$ is called a {\em mixed group identity} over $G$. 
When it is satisfied for all elements of $G$ we call it a {\em law with constants} of $G$. 
In this paper we study mixed identities/laws with constants in groups of homeomorphisms of some natural topological/metric spaces. 

Investigations of mixed identities were initiated by several Russian mathematicians in the 1970 - 80-th, see \cite{An1}, \cite{An2}, \cite{Gol} and \cite{To}.   
The topic was developed later in \cite{G}, \cite{NSt} and \cite{St}.  
Very recently groups without mixed identities ({\em MIF groups}) have become an attractive object of investigations under some geometric assumptions. 
In particular Hull and Osin have observed in \cite{HO} that acylindically hyperbolic groups with trivial finite radical are MIF. 
Papers \cite{HO} and \cite{Jacobson} contain many other  interesting examples in connection with some dynamical and geometric properties studied in group theory.  
We also mention an observation from \cite{HO} that MIF can be reformulated in some model-theoretic terms (and even in terms of algebraic geometry). 

Investigation of mixed identities is fruitful in the context of permutation groups. 
For example Theorem 1.6 of \cite{HO} states that in the class of highly transitive countable groups property MIF is opposite to embeddability of $\mathsf{Alt} (\mathbb{N})$ as a normal subgroup.  
It is shown in \cite{EGMM} that dense countable subgroups of automorphism/isometry groups of typical structures/spaces with some properties of universality, are MIF.    
On the other hand it is proved in \cite{Zar} and in \cite{SZ} that Thompson's group $F$ is not MIF.  

The latter fact together with the Abert's theorem that Thompson's group $F$ does not satisfy any identity have become the starting point of this paper. 
What are mixed identities which are not satisfied in Thompson's group $F$?  
Generalizing the Abert's property of separation we introduce hereditarily separating $G$-spaces and prove that under this condition the group $G$ does not satisfy so called {\em oscillating} identities with constants.  
Since the standard action of $F$ on $[0,1]$ is hereditarily separating we arrive at some partial answer to the above question. 
The same situation arises in the context of actions of branch groups on the boundary space. 
A short formulation of the main theorem is as follows (the exact formulation is given in Theorem 2.13 in Section 2 together with the corresponding  definitions). 
\bigskip 

\noindent 
{\bf Theorem 1.} 
{\em Let $G$ be a group acting on a perfect Hausdorff space $\mathcal{X}$ by homeomorphisms. 
Let $w(\bar{x})$ be a non-trivial reduced word from  
$\mathbb{F}_{n}\ast G$ such that $w(\bar{x})\notin G$. 
If $G$ hereditarily separates $\mathcal{X}$ and $w(\bar{x})$ 
has a conjugate in $\mathbb{F}_{n}\ast G$ which 
is explicitly oscillating, then the inequality $w(\bar{x})\neq 1$ has a solution in $G$. }

\bigskip 

\noindent 
Using this theorem we will also see that explicit oscillation gives a reasonable substitute for MIF, see Section 2.3. 
  
It is worth noting that the idea of oscillation which we introduce is straightforward: this is a condition on the family 
of the supports of the constants of the word. 
In particular our arguments are rather direct but technical.  
Examples which we give in the paper should justify this inconvenience: they show that the notion is really ubiquitous and unavoidable.

The major results of this paper appeared in some form in Section 2 of the PhD thesis of the second author (under supervision of the first one) in \cite{Zarphd}.  
In Sections 3 and 4 of \cite{Zarphd} they are applied to investigations of $F$-limit groups. 
These applications will form a separate article. 

We view the present article as an attempt to generalize the property MIF in the context of $G$-spaces.  
Theorem \ref{ab} gives a sufficient condition for the existence of a solution of a given inequality in a group having a hereditarily separating action. 
Techniques applied in the proof will be enriched in other parts of the paper in order to cover new cases under new circumstances. 
In particular, we develop the approach in Theorem \ref{uab}, which deals with a finite set of inequalities and uses weaker
assumptions. 
The case for which $G$ is equipped with a topology  is considered in Section 4. 

The terminology and notation used in the paper are standard.  
In the final part of the introduction we remind the reader of some information concerning Thompson's group $F$. 
In the main body of the paper it is assumed that the reader remembers it (in particular Fact below).  
We mainly follow \cite{CFP}.
\bigskip 

{\em Thompson's group} $F$ is the group given by the following infinite group presentation:
\[ 
\langle x_{0}, x_{1}, x_{2}, \ldots \Big| \, x_{j}x_{i}=x_{i}x_{j+1} \, i<j  \rangle . 
\] 
In fact $F$ is finitely presented:
\[ 
F = \langle x_{0}, x_{1}\ \Big|\ [x_{0}x_{1}^{-1}, x_{0}^{-i}x_{1}x_{0}^{i}] \, , \, i=1,2 \rangle .  
\] 
We will use the following geometric interpretation of $F$. 
Consider the set of all strictly increasing continuous piecewise-linear functions from the closed unit interval onto itself. 
Then the group $F$ is realized by the set of all such
functions, which are differentiable except at finitely many dyadic rational numbers and such that all slopes (deriviatives) are integer powers of 2. 
The corresponding group operation is just the composition. 
For the further reference it will be useful to give an explicit form of the generators $x_{n}$, for $n\geq 0$, in terms of piecewise-linear functions:
\[ 
x_{n}(t) = \left\{ \begin{array}{ll} t & \textrm{, $t\in [0,\frac{2^{n}-1}{2^{n}} ]$} \\
	\frac{1}{2}t + \frac{2^{n}-1}{2^{n+1}} & \textrm{, $t\in [\frac{2^{n}-1}{2^{n}}, \frac{2^{n+1}-1}
	{2^{n+1}} ]$} \\ t - \frac{1}{2^{n+2}} & \textrm{, $t\in [\frac{2^{n+1} -1}{2^{n+1}},
	\frac{2^{n+2}-1}{2^{n+2}}]$} \\	2t-1 & \textrm{, $t\in [\frac{2^{n+2}-1}{2^{n+2}},1] .$} \end{array}\right. 	 
\] 
For any dyadic subinterval $[a,b]\subset [0,1]$, let us consider the set, of elements in $F$, which are trivial on its complement, and denote it by $F _{[a,b]}$. 
We know that it forms a subgroup of $F$, which is isomorphic to the whole group. 
Let us denote its standard infinite set of generators by $x_{[a,b],0}, x_{[a,b],1}, x_{[a,b],2}, \ldots$, 
where for $n\geq 0$ we have:
\[ 
x_{[a,b],n}(t) = \left\{ \begin{array}{ll} t & \textrm{, $t\in [0,a+\frac{(2^{n}-1)(b-a)}{2^{n}}]$} \\
	\frac{1}{2}t + \frac{1}{2}(a+\frac{2^{n}-1}{2^{n}}) & \textrm{, $t\in [a+\frac{(2^{n}-1)(b-a)}{2^{n}},
			a+\frac{(2^{n+1}-1)(b-a)}{2^{n+1}}]$} \\
	 t - \frac{b-a}{2^{n+2}} & \textrm{, $t\in [a+\frac{(2^{n+1}-1)(b-a)}{2^{n+1}},
			a+\frac{(2^{n+2}-1)(b-a)}{2^{n+2}}]$} \\
	2t - b & \textrm{, $t\in [a+\frac{(2^{n+2}-1)(b-a)}{2^{n+2}},b]$} \\
	t & \textrm{, $t\in [b,1] .$} \end{array}\right. 
\] 
Moreover, if $\iota _{[a,b]}$ denotes the natural isomorphism between $F$ and $F_{[a,b]}$ sending $x_{n}$ to $x_{[a,b],n}$ for all $n\geq 0$, then for any $f\in F$ by $f_{[a,b]}$ we denote the element $\iota _{[a,b]} (f)\in F_{[a,b]}<F$. 

We will repeatedly use the following fact (see Lemmas 4.2 and 2.4 in 
\cite{CFP} and  \cite{KM} respectively). 

\bigskip 

\noindent 
{\bf Fact.} \label{CFP} 
{\em If $0=x_{0}<x_{1}<x_{2}<\ldots <x_{n}=1$ and 
$0=y_{0}<y_{1}<y_{2}<\ldots <y_{n}=1$ are partitions of $[0,1]$ consisting of dyadic rational numbers, then there exists $f\in F$ such that $f(x_{i})=y_{i}$ for $i=0,\ldots , n$. 
	
Furthermore, if $x_{i-1}=y_{i-1}$
and $x_{i}=y_{i}$ for some $i$ with $1\leq i\leq n$, then $f$ can be taken to be trivial on the interval $[x_{i-1},x_{i}]$.
}

\bigskip 

Many examples, which occur in this paper, can be easily exposed using the \emph{rectangle diagrams} introduced by W. Thurston. 
The paper \cite{Zarphd} contains the corresponding pictures.  

\section{Explicitly oscillating words and hereditarily separating actions}

In this section we study inequalities over groups of permutations which are similar to separating actions introduced by Abert in \cite{A}. 
 
Let $G$ be a group. 
Any inequality over $G$ can be considered as follows. 
Let $w(\bar{y})$ be a word over $G$ on $t$ variables 
$y_{1},\ldots , y_{t}$. 
We usually view it as a nontrivial element of $\mathbb{F} _{t}\ast G$. 
In order to study the inequality $w(\bar{y}) \not= 1$
we will assume that $w(\bar{y})$ is reduced in $\mathbb{F} _{t}\ast G$. 
If $w(\bar{y})\notin\mathbb{F} _{t}$, we usually assume that it is of the form 
\[ 
w(\bar{y})=\mathsf{u}_{n}v_{n}\mathsf{u}_{n-1}v_{n-1}\ldots \mathsf{u}_{1}v_{1}, \hspace{4cm}(1.1)
\]  
where $n\in\mathbb{N}$, $\mathsf{u}_{i}\in \mathbb{F}(\bar{y})$ 
and $v_{i}\in G\setminus\{ 1\}$ for each $i\leq n$. 
It is clear that any word with constants is conjugate to a word in this form. 

Our basic concern is existence of solutions of the inequalitiy 
$w(\bar{y}) \not=1$ in $G$. 
The easiest case appears in the following definition.   

\begin{definition} \label{prod-const} 
If in the form {\em (1.1)}  $v_{n}\cdot \ldots \cdot v_{1}\neq 1$ then we say that the word $w(\bar{y})$ has non-trivial product of constants  
(in $\mathsf{supp}(v_{n}\cdot \ldots \cdot v_{1})\subseteq X$).
\end{definition} 
Then it is clear that the tuple of units $\bar{1}$ solves the inequalitiy  $w(\bar{y}) \not=1$. 

\begin{remark} 
{\em The subgroup of $G\ast\mathbb{F}_{t}$	consisting of words having trivial product of constants can be viewed as follows. 
The $t$-tuple of units $\bar{1}\in G^t$ is the variety of  
the system of equations $S_{triv}$: $y_{i}=1$, $1\leq i\leq t$.
Then according to Definition 2 of \cite{KhM} the set of all words $w(\bar{y}) \in \mathbb{F} _{t}\ast G$ with $w(\bar{1})=1$ is the radical of the system $S_{triv}$. }
\end{remark}
In order to study more general/interesting cases we will concentrate on permutation groups.  
In this class of groups we will introduce oscillating words. 
It will turn out that they generalize words with non-trivial product of constants.

\subsection{Explicitly oscillating words} 
Let $G$ be a permutation group on $X$. 
We distinguish some specific types of words over $G$ with
respect to the action on $X$. 
Let $w(\bar{y})$ be in the form (1.1). 
Define: 
\[ 
O_{w}:=\bigcap _{i=0} ^{n-1} v_{0}^{-1}v_{1}^{-1}\ldots v_{i}^{-1}( \mathsf{supp}(v_{i+1})) ,
\] 
where $v_{0}=\mathsf{id}$. 
If $w(\bar{y})\in\mathbb{F}_{t}$ then let 
$O_{w}:=X\setminus \mathsf{Fix}(G)$.

\begin{definition} \label{osc} 
Let $V\subseteq X$. 
We say that the word $w(\bar{y})\in\mathbb{F} _{t}\ast G$ 
is  explicitly oscillating in $V$ if $w(\bar{y})$ is non-trivial and $V\cap O_{w}\neq\emptyset$. 

When $V=X$ we just say that $w(\bar{y})$ is explicitly oscillating.
\end{definition} 
 
Observe that any non-trivial initial/final segment of an explicitly oscillating word is also explicitly oscillating. 
Note also that if all $v_{i}$ are taken from the same cyclic subgroup of $G$ of prime (or infinite) order, then $w(\bar{y})$ is explicitly oscillating. 
It is also clear that conjugating an explicitly oscillating word $w(\bar{y})$ with $n>1$ by an element of $G$ with support disjoint from $\bigcup _{i=1}^{n} \mathsf{supp}(v_{i})$ we obtain a word which is not explicitly oscillating. 

In the situation where $G$ acts on a Hausdorff topological space $\mathcal{X}$ by homeomorphisms the set $O_{w}$ is open. 
In the main results of the paper we usually assume that $\mathcal{X}$ is perfect. 
This rules out explicitly oscillating words of the following example. 

\begin{example} \label{fin} 
{\em Consider the group of finitary permutations 
$G=\mathsf{S_{fin}}(\mathbb{N})$. 
Let $w(\bar{y})\notin\mathbb{F}_{t}$ and $w(\bar{y})$ be in the form (1.1) with non-trivial $v_i$.  
Then $O_w$ is finite. }
\end{example}

\begin{remark}
{\em Nevertheless we will show in several places of the paper how our results can be adapted to the discrete case, see for example Theorem \ref{gab}. } 
\end{remark}

The case of Thompson's group $F$ is fundamental for us. 
It will serve as the main source of illustrations. 
They are usually based on some computations which are left to the reader. 
Details of these computations and the corresponding pictures are given in the PhD thesis of the second author, see \cite{Zarphd}.    

\begin{example} \label{1}
\emph{ Consider Thompson's group $F$ with its standard action on $[0,1]$. 
Let $w_{1}(y) =yx_{1}y^{-1}x_{2}y^{2}x_{1}^{-1}$, where $y$ is a variable.
In the notation from the definitions we have $v_{1}=x_{1}^{-1}$, $v_{2}=x_{2}$, $v_{3}=x_{1}$ and hence:
\[ 
O_{w_{1}}=x_{1}x_{2}^{-1}\Big( (\frac{1}{2},1 )\Big)\cap x_{1}\Big( (\frac{3}{4},1) \Big)\cap (\frac{1}{2},1) = (\frac{5}{8},1) . 
\] 
Thus we see that $w_{1}(y)$ is explicitly oscillating 
(in $(\frac{5}{8},1)$).}
\end{example}

\noindent 
{\bf Notation.} 
We introduce the following notation. 
For an element $v\in G$ let $v^1 = v$ and $v^0 = 1$.   
For every $w(\bar{y})$ given in the form 
$w(\bar{y}) =\mathsf{u}_{n}v_{n}\mathsf{u}_{n-1}v_{n-1}\ldots \mathsf{u}_{1}v_{1}$ 
as in (1.1) and for every set $A\subseteq X$ let 
\[ 
\mathcal{V}_{w}(A) = \{  v^{\varepsilon_j}_{j}\cdot \ldots \cdot v^{\varepsilon_1}_{1}(A)\, | \, (\varepsilon_1 ,\ldots , \varepsilon_j ) \in \{ 0,1 \}^j \, , \, 1 \le j \le n \} , 
\] 
\[ 
\mathcal{V}^{-1}_{w}(A) = \{  v^{\varepsilon_1}_{1}\cdot \ldots \cdot v^{\varepsilon_j}_{j}(A)\, | \, (\varepsilon_1 ,\ldots , \varepsilon_j ) \in \{ 0,-1 \}^j \, , \, 1 \le j \le n \} ,  
\] 
and let  
\[ 
\mathcal{V}^{>0}_{w}(A)= \{ v_{j}\cdot \ldots \cdot v_{1}(A)\, |\, 1\le j \le n \}.
\]   
When $w(\bar{y}) \in \mathbb{F}_{t}$ let $\mathcal{V}^{>0}_{w}(A)=\{ A\}$. 

\subsection{Hereditarily separating actions by homeomorphisms} 

\bigskip 

\begin{definition}
	Let $G$ be a permutation group acting on an infinite set $X$. \begin{itemize} 
\item We say that $G$ separates $X$ if, for every finite subset $Y\subset X$, the pointwise stabilizer 
$\mathsf{stab}_{G}(Y)$ does not stabilize any point outside $Y$. 
\item Assume that $\mathcal{X}$ is a Hausdorff topological space, $G$ consists of homeomorphisms of $\mathcal{X}$ and the set of fixed points, $\mathsf{Fix}(G)$, is finite. 
We say that $G$ hereditarily separates $\mathcal{X}$ if, for every  open and infinite subset $Z\subseteq \mathcal{X}$ and for every finite subset $Y\subset Z$, the subgroup 
$\mathsf{stab}_{G}((\mathcal{X}\setminus Z)\cup Y)$ 
does not stabilize any point from
	$Z\setminus (Y\cup \mathsf{Fix}(G))$.
\end{itemize} 
\end{definition}

\begin{remark} \label{rema}
{\em Separating actions were introduced in \cite{A}. 
We mention the following observation from that paper. 	 
For every separating action of a group $G$ on $X$ and every finite $Y\subset X$ the orbits of the action of the pointwise stabilizer of $Y$ in $X \setminus Y$ are infinite. 
As a corollary we see that }  
\begin{itemize} 
\item {\em for a hereditarily separating action of $G$ on $\mathcal{X}$ and an open and infinite subset $Z\subseteq \mathcal{X}$, the action of the stabilizer $\mathsf{stab}_{G}(\mathcal{X}\setminus Z)$ on $Z\setminus \mathsf{Fix}(G)$ has only infinite orbits.  }  
\end{itemize} 
\end{remark}

\begin{example} \label{alt} 
{\em Consider $\mathbb{N}$ as the discrete space. 
Then the symmetric group $\mathsf{Sym}(\mathbb{N})$, the group of finitary permutations $\mathsf{S_{fin}}(\mathbb{N})$ and the alternating group 
$\mathsf{Alt}(\mathbb{N}) \le \mathsf{S_{fin}}(\mathbb{N})$ hereditarily separate $\mathbb{N}$. }
\end{example}

\begin{example} \label{f}
\emph{ Observe that} Thompson's group $F$ is hereditarily separating with respect to its standard action on $[0,1]$. 
{\em Indeed suppose that $Z$ is an open subset of 
$[0,1]$, $Y:=\{ t_{1}, t_{2}, \ldots t_{s}\}\subset Z$ 
and $t\in Z\setminus (Y\cup \{ 0, 1\})$.
There is some non-trivial dyadic segment $[p,q]\subseteq Z$ containing $t$ such that $Y	\cap [p,q]=\emptyset$. 
Obviously 
$x_{[p,q],0}\in \mathsf{stab}_{F} (([0,1]\setminus Z)\cup Y)
	\setminus \mathsf{stab}_{F} (\{ t\} )$.}
\end{example}

\begin{example} \label{grig}  
The action of any weakly branch group on the boundary space of the corresponding infinite rooted tree is also hereditarily separating. 	
{\em To define this class we will follow \cite{Gr} and \cite{Gri}. 
We consider a group $G$, which isometrically acts on some rooted tree $T$. 
The vertices in $T$, which are at the same distance from the root are said to be at the same \emph{level}. 
We say that the action of the group $G$ on $T$ is} 
spherically transitive 
{\em if $G$ acts transitively on each level of $T$. 
For any vertex $t\in T$ we define its} 
rigid stabilizer 
{\em to be the set of all elements from $G$, which stabilize  
$T\setminus T_t$ pointwise, where $T_t$ is the natural subtree hanging from $t$. 
A group $G$ is called} 
weakly branch  
{\em if it acts spherically transitively on some rooted tree $T$ so that the rigid stabilizer of every vertex is non-trivial.
The} 
boundary 
{\em of a tree $T$, denoted by $\partial T$, consists of the infinite branches starting at the root. 
There is a topology on $\partial T$ where the base of open sets is determined by subtrees $T_t$. 
A weakly branch group obviously acts on the boundary $\partial T$ by homeomorphisms (in fact by isometries with respect to the natural metric). 

To see the statement of this example we use an argument from \cite{A}.  
Fix a weakly branch group $G$ and the corresponding tree $T$. 
Let $\mathcal{X} :=\partial T$. 
To see that $G$ hereditarily separates $\mathcal{X}$, suppose that $Z$ is an open	subset of	$\mathcal{X}$, 
$Y:= \{ x_{1}, x_{2}, \ldots x_{\ell}\} \subset Z$ and 
$x\in Z\setminus Y$.
Wlog assume that $t$ is the vertex of the infinite ray $x$ such that if any $y\in \mathcal{X}$ contains $t$ then $y\in Z$. 
Let $k$ be the level of $t$. 
Now choose a level $k'\geq k$ such that the vertices in the infinite rays $x_{1},\ldots , x_{\ell}, x$ at the $k'$-th level 
are all distinct. 
Let $t_{0}$ be the vertex of $x$ at level $k'$. 
Let $S$ be the stabilizer of $t_{0}$ in $G$ and let $R$ be the rigid vertex stabilizer of $t_{0}$ in $G$. 
Then $S$ acts spherically transitively on the infinite subtree $T_{t_{0}}$ rooted at $t_{0}$.} 
{\em Since $R$ is a non-trivial normal subgroup of $S$ we see that it cannot stabilize any infinite ray going through $t_{0}$. 
In particular, there exists $g\in R$ such that $g(x)\neq x$. 
On	the other hand, $g$ stabilizes every ray not going through $t_{0}$.
It follows that $g\in \mathsf{stab} _{G}((\mathcal{X}\setminus Z)\cup Y)$.
This proves the statement of Example \ref{grig}. }
\end{example}

Before the formulation of Theorem \ref{ab} we introduce a notion which will be central in the proof. 
Let $G$ be a permutation group on $X$ and  
$w(\bar{y})$ be a word over $G$ on $t$ variables $y_{1},\ldots , y_{t}$. 
Assume that $w(\bar{y})$ is reduced in $\mathbb{F} _{t}\ast G$,  
$w(\bar{y})\notin\mathbb{F} _{t}$ and  
$w(\bar{y})$ is written in the following form:
\[ 
w(\bar{y}) =u_{n,\ell_{n}}\ldots u_{n,1}v_{n}\ldots u_{2,\ell_{2}}\ldots u_{2,1}v_{2}u_{1,\ell_{1}}\ldots  u_{1,1}v_{1}, \hspace{2cm}(1.2)
\]
where 
$u_{j,i_{j}}\in\{ y_{1}^{\pm 1},\ldots , y_{t}^{\pm 1}\}$, $1\leq i_{j}\leq \ell_{j}$, and 
$v_{j}\in G \setminus\{ 1\}$, $1\leq j\leq n$. 
Let $L_{j}:=\sum_{i=1}^{j} \ell_{i}$, $1\leq j\leq n$. 

From now on we will say that $L_n$ is the {\em length} of $w(\bar{y})$. 
It will be denoted by $|w(\bar{y})|$. 
This terminology slightly disagrees with the usual one (where parameters are taken into account), but it will be helpful in our induction arguments. 

For any $1\leq r\leq L_{n}$ let  
\[ 
(w)_{r}(\bar{y}) =u_{d,s}\ldots u_{d,1}v_{d}\ldots u_{2,\ell_{2}}\ldots u_{2,1}v_{2}u_{1,\ell_{1}}\ldots  u_{1,1}v_{1}  
\]
be the final segment of $w(\bar{y})$ of $r$ occurances of letters from $\{ y_{1}^{\pm 1},\ldots , y_{t}^{\pm 1}\}$, 
i.e. $r=L_{d-1}+s$ and $1\leq s\leq \ell_{d}$. 
In this case we also define  
\[ 
[w]_{r}(\bar{y}) =w(\bar{y}) ((w)_{r}(\bar{y}))^{-1}  
\]
to be the corresponding initial segment. 

Let $\bar{g} = (g_{1},\ldots , g_{t})$ be a tuple from $G$. 
By $(w)_{r}(\bar{g})$ we denote the value of $(w)_{r}(\bar{y})$ in $G$ under the substitution $y_{i}=g_{i}$, $1\leq i\leq t$. 
To simplify notation let also $(w)_{0}(\bar{y})=1\in G$. 

Let $p\in X$. 
Define $p_{r,\bar{g}}=(w)_{r}(\bar{g})(p)$ for all $r$, $1\leq r\leq L_{n}$. 

\begin{definition} \label{dist} 
We say that $\bar{g}$ is distinctive for $w(\bar{y})$ and $p$, if all the points
\[ 
p=p_{0,\bar{g}}, v_{1}(p_{0,\bar{g}}),\ldots , p_{l_{1},\bar{g}}, v_{2}(p_{l_{1},\bar{g}}),\ldots , p_{L_n , \bar{g}}
\] 
are pairwise distinct. 
\end{definition} 
This generalizes the corresponding definition of \cite{A} given on p.528 for words $w(\bar{y}) \in \mathbb{F}_t$. 
In the latter case we just omit the members  
$v_{j}(p_{l_{j-1},\bar{g}})$ from the sequence. 
When Definition \ref{dist} holds we obviously have 
$w(\bar{g})(p)\neq p$, i.e. $w(\bar{g}) \not=1$. 
The following statement is related to Theorem 1.1 from \cite{A}.

\begin{theorem} \label{ab}
Let $G$ be a group acting on a perfect Hausdorff space $\mathcal{X}$ by homeomorphisms. 
Let $w(\bar{y})$ be a nontrivial word over $G$ on $t$ variables, $y_{1},\ldots , y_{t}$, which is reduced in 
$\mathbb{F}_{t}\ast G$ and non-constant 
(i.e. $w(\bar{y})\notin G$). 
If $G$ hereditarily separates $\mathcal{X}$ and $w(\bar{y})$ 
has a conjugate in $\mathbb{F}_{t}\ast G$ which 
is explicitly oscillating, then the inequality $w(\bar{y})\neq 1$ has a solution in $G$. 

Moreover, assuming additionally that 
$w(\bar{y})\in\mathbb{F} _{t}$ or $w(\bar{y})$ is given in the form $(1.1)$ and is explicitly oscillating,  for every open set $O'\subseteq O_{w}$, there is a distinctive tuple $\bar{g}=(g_{1},\ldots , g_{t})\in G^t$ with respect to $w(\bar{y})$ and some $p\in O'$ such that for all $i$ with $1\leq i\leq t$, 
\begin{itemize} 
\item $\mathsf{supp}(g_{i})\subseteq (\bigcup \mathcal{V}^{>0}_{w}(O'))\setminus (\bigcup \mathcal{V}_w (\overline{O'}\setminus O'))$, 
\item each member of $\mathcal{V} _{w}(O')$ is $g_i$-invariant.  
\end{itemize} 
\end{theorem}

\emph{Proof.} 
Using conjugation (if necessary) we assume that if 
$w(\bar{y})\notin\mathbb{F} _{t}$ then $w(\bar{y})$ is 
an explicitly oscillating word written in the form (1.2). 
We preserve the notation given before the formulation. 
Note that since $\mathcal{X}$ is perfect, each open subset is infinite. 

It is clear that we only need to prove the second part of the statement of the theorem.  
Fix $O'$ as in the formulation and any 
$p\in O'\setminus (\bigcup \mathcal{V}^{-1}_w (\bigcup \mathcal{V}_w (\overline{O'}\setminus O')))$. 
The existence of such $p$ follows from the fact that the set 
$\overline{O'}\setminus O'$
is nowhere dense in $\mathcal{X}$ and the action of $G$ is continuous. 
The proof of the theorem is by induction. 
At $k$-th step (where $k\le L_n$) we will show that: 
\begin{itemize} 
\item There is a tuple 
$\bar{g}=(g_{1},\ldots , g_{t})\in G$ such that $\bar{g}$ is distinctive for $(w)_{k}$ and $p$. 
\item In the condition above we can choose $\bar{g}$ so that for all $i$ with $1\leq i\leq t$, 
$\mathsf{supp}(g_{i})\subseteq (\bigcup \mathcal{V}^{>0}_w (O'))\setminus (\bigcup \mathcal{V}^{-1}_{[w]_{k}} (\bigcup \mathcal{V}_w (\overline{O'}\setminus O')))$ and each member of $\mathcal{V}_{w}(O')$ is $g_i$-invariant. 
\end{itemize} 
We will use the following claim.

\noindent  
 \textbf{Claim M.} 
{\em For every $r\leq n$, every 
$(\varepsilon_1 \ldots \varepsilon_r)\in \{ 0,1 \}^r$  and every  
\[ 
q\in v^{\varepsilon_r}_{r}\cdot \ldots \cdot v^{\varepsilon_1}_{1}(O')\setminus\bigcup  \mathcal{V}_{w}(\overline{O'}\setminus O') 
\] 
there is a neighbourhood $O\subseteq v^{\varepsilon_r}_{r}\cdot \ldots \cdot v^{\varepsilon_1}_{1}(O')$ of $q$ such that: }
\[  
 \forall V\in \mathcal{V}_{w}(O') 
\bigg( O\cap V\neq\emptyset\ \Rightarrow\ O\subseteq V\bigg) . \hspace{3cm}(\dag )
\]
Indeed, let $q$ and $O'$ be as in the formulation of the claim. 
Let us enumerate the elements of $\mathcal{V}_{w}(O')$ by $V_m$ with $m>0$ and construct inductively open sets $O_m$ corresponding to them. 
We start with 
$O_{0}=v^{\varepsilon_r}_{r}\cdot \ldots \cdot v^{\varepsilon_1}_{1}(O')$. 
If at the $m$-th step $q\in O_{m-1}\cap V_m$ then we define $O_{m} =O_{m-1}\cap V_m$. 
Now suppose that $q\notin O_{m-1}\cap V_m$. 
Since 
$q\notin \bigcup  \mathcal{V}_{w}(\overline{O'}\setminus O')$,
we find some open  
$O_{m}\subseteq O_{m-1}\setminus V_m$ containing $q$. 
After meeting all members of $\mathcal{V}_{w}(O')$  
we obtain the final $O$ which is non-empty and open. 
If $w(\bar{y})\in\mathbb{F}_{t}$ then we define $O=O'$, which finishes the proof of the claim. 

\bigskip 

If $k=1$, we may assume that $(w)_{1}$ is of the form 
$y_{j}^{\pm 1}v_{1}$ for some $1\leq j\leq t$. 
When $w(\bar{y})\in\mathbb{F} _{t}$ we replace $v_{1}$ by $\mathsf{id}$ and follow the argument below. 
Note that 
$v_{1}(p)\notin \bigcup \mathcal{V}^{-1}_{[w]_{1}}(\bigcup  \mathcal{V}_{w}(\overline{O'}\setminus O'))$. 
Applying the claim above we find in $v_{1}(O')$ a neighborhood 
$O$ of $v_{1}(p)$ which satisfies $(\dag )$.
According to the assumptions above, 
when $w(\bar{y})\not\in\mathbb{F} _{t}$ then 
$p'\neq v_{1}(p')$ for all  $p'\in O'$. 
Thus for a non-trivial $v_1$ 
the inequality $p\neq v_{1}(p)$ is satisfied. 

Wlog suppose $(w)_1 =y_{j}v_{1}$. 
Since $G$ hereditarily separates $\mathcal{X}$, the orbit of $v_{1}(p)$ with respect to $\mathsf{stab}_{G}((\mathcal{X}\setminus O)\cup \bigcup \mathcal{V}^{-1}_{[w]_1} (\bigcup \mathcal{V}_w (\overline{O'}\setminus O')))$ is infinite. 
Thus we can choose $f\in \mathsf{stab}_{G}((\mathcal{X}\setminus O)\cup \bigcup \mathcal{V}^{-1}_{[w]_1} (\bigcup \mathcal{V}_w (\overline{O'}\setminus O')))$ such that $f(v_{1}(p))\notin\{ p, v_{1}(p)\}$. 
Defining $g_{j}:=f$ and extending it by any $(t-1)$-tuple from 
$\mathsf{stab}_{G}((\mathcal{X}\setminus O)\cup \bigcup \mathcal{V}^{-1}_{[w]_1} (\bigcup \mathcal{V}_w (\overline{O'}\setminus O')))$ we obtain a $t$-tuple $\bar{g}$ distinctive for $p$ and $(w)_{1}$. 

It is easy to see that for all $i\leq t$,
\[ 
\mathsf{supp}(g_{i})\subseteq (\bigcup \mathcal{V}^{>0}_{w}(O'))\setminus (\bigcup \mathcal{V}^{-1}_{[w]_1} (\bigcup \mathcal{V}_w (\overline{O'}\setminus O'))). 
\] 
Now fix $i\leq t$ and $V\in \mathcal{V}_{w}(O')$. 
If $O\cap V=\emptyset$, then obviously  $g_{i}(V)=V$. 
On the other hand, if $O\cap V\neq\emptyset$, then it follows from the construction of $O$ that $O\subseteq V$. 
Thus $\mathsf{supp}(g_{i})$ is also a subset of  $V$ and the condition  $g_{i}(V)=V$ is satisfied. 

Let $k>1$. 
Assume that for
\[ 
(w)_{k-1}=u_{d,s}\ldots u_{d,1}v_{d}\ldots u_{2, \ell_2}\ldots u_{2,1}v_{2}u_{1,\ell_{1}}\ldots u_{1,1}v_{1} ,
\]
where $k-1=L_{d-1}+s$, we can find a tuple 
$\bar{g}\in G$ such that the induction hypothesis is satisfied. 
According to the form of $(w)_{k}$ we consider two cases. 

\noindent
{\em Case 1.} $(w)_{k}=u_{d,s+1}u_{d,s}\ldots u_{d,1}v_{d}\ldots u_{2,1}v_{2}u_{1,\ell_{1}}\ldots u_{1,1}v_{1}$, where $k-1=L_{d-1}+s$, $s\geq 1$.\parskip0pt

 If
\[ 
p_{k,\bar{g}}\notin\Big\{ p_{i,\bar{g}}\ \Big|\ 0\leq i\leq k-1\Big\}\cup\Big\{ v_{1}(p_{0,\bar{g}}),\ldots
	, v_{d}(p_{L_{d-1},\bar{g}})\Big\} , 
\] 
then we have found an acceptable tuple $\bar{g}$. 
Let us assume that $p_{k,\bar{g}}=p_{m,\bar{g}}$ for some
$0\leq m<k$ or $p_{k,\bar{g}}=v_{m+1}(p_{L_{m},\bar{g}})$ for some $0\leq m<d-1$. 

 Let $y_{j}^{\pm 1}$ be the first letter of $(w)_{k}$. 
 Replacing $y_{j}$ by $y_{j}^{-1}$ and $g_{j}$ by
$g_{j}^{-1}$ if necessary, we may assume that $u_{d,s+1}=y_{j}$. 
Put
\[ 
Y =\Big\{ p_{i,\bar{g}}\ \Big|\ 0\leq i\leq k-2\Big\}\cup\Big\{ v_{1}(p_{0,\bar{g}}),\ldots , v_{d}(p_{L_{d-1},
	\bar{g}})\Big\} 
\] 
for $\bar{g}$ chosen at the $(k-1)$-th step of induction. 
Since by the induction hypothesis $p\in O'$ and 
$g_{s}^{\pm 1}(v_{r}\ldots v_{1}(O'))=v_{r}\ldots v_{1}(O')$ for all $r\leq k$ and all $ s\leq t$, we see that 
$p_{k-1,\bar{g}}\in v_{d}\ldots v_{1}(O')$. 
We also know that 
$p_{k-1,\bar{g}}\notin\bigcup \mathcal{V}^{-1}_{[w]_{k-1}} (\bigcup \mathcal{V}_w (\overline{O'}\setminus O'))$. 
As above we choose a neighborhood 
$O\subseteq v_{d}\ldots v_{1}(O')$ of the point 
$p_{k-1,\bar{g}}\in v_{d}\ldots v_{1}(O')$ satisfying condition $(\dag )$. 
Since the action of $G$ is hereditarily separating, 
the 
$\mathsf{stab}_{G}((\mathcal{X}\setminus O)\cup \bigcup \mathcal{V}^{-1}_{[w]_{k-1}} (\bigcup \mathcal{V}_w (\overline{O'}\setminus O')))$-orbit of $p_{k-1,\bar{g}}$ is infinite. 
Let 
\[ 
Z =\Big\{ g_{j}^{-1}(p_{i,\bar{g}})\ \Big| \ 0\leq i\leq k-1\Big\}\cup\Big\{ g_{j}^{-1}(v_{i+1}(p_{L_{i},\bar{g}}))\ \Big|\ 0\leq i\leq d-1\Big\} . 
\] 
Since $Z$ is finite, there exists $f\in \mathsf{stab}_{G}((\mathcal{X}\setminus O)\cup \bigcup \mathcal{V}^{-1}_{[w]_{k-1}} (\bigcup \mathcal{V}_w (\overline{O'}\setminus O'))\cup Y)$ taking $p_{k-1,\bar{g}}$ outside $Z$. 
Replacing $g_{j}$ by $g_{j}f$ we obtain a corrected tuple $\bar{g}$. 
Since the element $f$ has been chosen from the stabilizer $\mathsf{stab}_{G}(Y)$, the points
\[ 
p_{0,\bar{g}},v_{1}(p_{0,\bar{g}}), p_{1,\bar{g}},\ldots , v_{d}(p_{L_{d-1},\bar{g}}), p_{L_{d-1}, \bar{g}},\ldots , p_{k-1,\bar{g}} 
\]
are the same as before. 
On the other hand $p_{k,\bar{g}}$ is distinct from all elements
\[ 
p_{0,\bar{g}},v_{1}(p_{0,\bar{g}}), p_{1,\bar{g}},\ldots , v_{d}(p_{L_{d-1},\bar{g}}), p_{L_{d-1}, \bar{g}},\ldots, p_{k-1,\bar{g}} . 
\] 
Since 
\[ 
f\in \mathsf{stab}_{G}((X\setminus O)\cup \bigcup \mathcal{V}^{-1}_{[w]_{k-1}} (\bigcup \mathcal{V}_w (\overline{O'}\setminus O'))\cup Y), 
\] 
we see 
$\mathsf{supp}(f)\subseteq\bigcup \mathcal{V}^{>0}_{w}(O')\setminus (\bigcup \mathcal{V}^{-1}_{[w]_{k-1}} (\bigcup \mathcal{V}_w (\overline{O'}\setminus O')))$. 
 This together with the induction hypothesis implies that 
\[ 
\mathsf{supp}(g_{j})\subseteq\bigcup\mathcal{V}^{>0}_{w}(O')\setminus (\bigcup \mathcal{V}^{-1}_{[w]_{k-1}} (\bigcup \mathcal{V}_w (\overline{O'}\setminus O')))
\]  
for the corrected $g_{j}$. 

Let $V\in \mathcal{V}_{w}(O')$. 
If $O\cap V=\emptyset$, then by the choice of $f$ and induction we have $g_{j}(V)=V$. 
If $O\cap V\neq\emptyset$, then $O\subseteq V$. 
Thus by the choice of $f$ and induction, $f$ and the original $g_{j}$ (defined at Step $k-1$) stabilize $V$ setwise. 
Thus $g_{j}(V)=V$. 
On the other hand all elements $g_{s}$ for $s\neq j$ have not been changed and thus automatically satisfy the required conditions. 
This finishes the proof of Case 1. 

\noindent 
{\em Case 2.} 
$(w)_{k}=u_{d+1,1}v_{d+1}u_{d,s}\ldots u_{d,1}v_{d}\ldots u_{2,1}v_{2}u_{1, \ell_{1}}\ldots u_{1,1}v_{1}$, 
where $k=L_{d}+1$. 

If
\[ 
p_{k,\bar{g}} \not\in \Big\{ p_{i,\bar{g}}\ \Big|\ 0\leq i\leq k-1\Big\}\cup \Big\{ v_{1}(p_{0,\bar{g}}),\ldots , v_{d+1}(p_{L_{d},\bar{g}})\Big\} , 
\]
and
\[ 
v_{d+1}(p_{k-1,\bar{g}})\notin\Big\{ p_{i,\bar{g}}\ \Big|\ 0\leq i\leq k-1\Big\}\cup\Big\{ v_{1}(p_{0,\bar{g}}),\ldots , v_{d}(p_{L_{d-1},\bar{g}})\Big\} , 
\]
then we have found an acceptable tuple $\bar{g}$. 

Assume the contrary. 
Suppose that $y_{j}^{\pm 1}$ is the first letter of $(w)_{k}$ 
and as in Case 1 we only consider the possibility 
$u_{d+1,1}=y_{j}$. 
Let $u_{d,s}=y_{j'}^{\pm 1}$. 
Then let
\[ 
Y' =\Big\{ p_{i,\bar{g}}\ \Big|\ 0\leq i\leq k-1\Big\}\cup\Big\{ v_{1}(p_{0,\bar{g}}),\ldots , v_{d}(p_{L_{d-1},\bar{g}})\Big\}. 
\] 
Assume that $v_{d+1}(p_{k-1,\bar{g}})\in Y'$. 
By the induction hypothesis, 
$p_{k-1,\bar{g}}\in v_{d}\ldots v_{1}(O')$. 
By Claim M we find some neighborhood 
$O\subseteq v_{d}\ldots v_{1}(O')$ of $p_{k-1,\bar{g}}$ 
satisfying $(\dag )$. 
Since the action of $G$ is hereditarily separating, the orbit of $p_{k-1,\bar{g}}$ with respect to 
\[ 
\mathsf{stab}_{G}((\mathcal{X}\setminus O)\cup \bigcup \mathcal{V}^{-1}_{[w]_{k-1}} (\bigcup \mathcal{V}_w (\overline{O'}\setminus O'))\cup (Y'\setminus\{ p_{k-1,\bar{g}}\} )) 
\] 
is infinite. 
We replace $g_{j'}$ by some $f'g_{j'}$ (or $g_{j'}f'$ in the case $u_{d,s}=y_{j'}^{-1}$), where
\[ 
f'\in \mathsf{stab}_{G}\Big( (\mathcal{X}\setminus O)\cup \bigcup \mathcal{V}^{-1}_{[w]_{k-1}} (\bigcup \mathcal{V}_w (\overline{O'}\setminus O'))\cup (Y'\setminus\{ p_{k-1,\bar{g}}\})\Big)
\] 
and $f'$ takes $p_{k-1,\bar{g}}$ outside the finite set 
$Y'\cup v_{d+1}^{-1}(Y')$. 
By the choice of $O'$ and $O$,    
\[ 
O\subseteq v_{d}\ldots v_{1}(O')\subseteq \mathsf{supp}(v_{d+1}) , 
\] 
i.e. the corrected $p_{k-1,\bar{g}}$ is not fixed by $v_{d+1}$. 
Thus the corrected $v_{d+1}(p_{k-1,\bar{g}})$ surely omits the corrected $Y'$. 

Furthermore, by the choice of $f'$ we still have 
\[ 
\mathsf{supp}(g_{j'})\subseteq\bigcup \mathcal{V}^{>0}_{w}(O')
\setminus (\bigcup \mathcal{V}^{-1}_{[w]_{k-1}} (\bigcup \mathcal{V}_w (\overline{O'}\setminus O'))).
\] 
Similarly as in Case 1 we see that the condition 
$\mathsf{supp}(f')\subseteq O$ implies that 
$g_{j'}^{\pm 1}(V)=V$ for any $V\in \mathcal{V}_{w}(O')$. 

So it remains to consider the case when 
$v_{d+1}(p_{k-1,\bar{g}})\notin Y'$, but either 
$p_{k,\bar{g}}=p_{j,\bar{g}}$ for some 
$0\leq j<k$ or $p_{k,\bar{g}}=v_{j+1}(p_{L_{j},\bar{g}})$ for some $0\leq j\leq d$. 
Let 
\[ 
Y =\Big\{ p_{i,\bar{g}}\ \Big|\ 0\leq i\leq k-1\Big\}\cup\Big\{ v_{1}(p_{0,\bar{g}}),\ldots , v_{d}(p_{L_{d-1},\bar{g}})\Big\}
\] 
(i.e. redefine $Y'$) and
\[ 
Z =\Big\{g_{j}^{-1}(p_{i,\bar{g}})\ \Big| \ 0\leq i\leq k-1\Big\}\cup\Big\{ g_{j}^{-1}(v_{i+1}(p_{L_{i}, \bar{g}}))\ \Big|\ 0\leq	 i\leq d\Big\} . 
\] 
Now observe that 
\[ 
v_{d+1}(p_{k-1,\bar{g}})\in v_{d+1}\ldots v_{1}(O')\setminus ( \bigcup \mathcal{V}^{-1}_{[w]_{k}} (\bigcup \mathcal{V}_w (\overline{O'}\setminus O'))\cup Y). 
\] 
Choose the neighborhood $O\subseteq v_{d+1}\ldots v_{1}(O')$ of the point $v_{d+1}(p_{k-1,\bar{g}})$ satisfying $(\dag )$. 
Then there exists $f\in \mathsf{stab}_{G}((\mathcal{X}\setminus O)\cup \bigcup \mathcal{V}^{-1}_{[w]_{k}} (\bigcup \mathcal{V}_w (\overline{O'}\setminus O'))\cup Y)$ taking $v_{d+1}(p_{k-1,\bar{g}})$ outside $Z$. 
Replacing $g_{j}$ by $g_{j}f$ we finish the proof exactly as in Case 1. $\square$

\bigskip 

\begin{remark} \label{diameter1} 
{\em Under the circumstances of Theorem \ref{ab} additionally assume that $(\mathcal{X},d)$ is a metric space. 
Then the conclusion of Theorem \ref{ab} may be extended by the following statement:  
\begin{itemize} 
\item for any $\varepsilon$ the tuple $\bar{g}$ can be chosen so that it additionally satisfies the inequality 
$d(x,g_i(x))\le \varepsilon$ for all $i\le t$ and $x\in \mathcal{X}$.  
\end{itemize}  
To see this it is enough to control how $\bar{g}$ is corrected at each step of the inductive procedure of the proof of Theorem \ref{ab}. 
If $|w(\bar{y})|=L_n$ (as above) then at $k$-th step of the induction we choose the corresponding open set $O$ with the additional property that its diameter is less than  $\frac{\varepsilon}{n}$. 
Assume that the homeomorphism $g_i$ which is corrected at this step already satisfies $d(x,g_i(x))\le \frac{k - 1}{n}\varepsilon$ for all $x\in \mathcal{X}$. 
Since the result of the correction, say $g'_i$, differs from $g_i$  by some $f'$ with the support from $O$ we have that 
$d(x , g'_i (x))\le \frac{k}{n} \varepsilon$ for all $x\in X$. } 
\end{remark}

We also have a topology-free version of Theorem \ref{ab}, which generalizes Theorem 1.1 from \cite{A}. 
For any $A\subseteq X$ we denote  $A^{1} = A$ and 
$A^0 = X\setminus A$. 
Fix any $w(\bar{y})$ such that either 
$w(\bar{y})\in\mathbb{F} _{t}$ or 
$w(\bar{y}) =\mathsf{u}_{n}v_{n}\ldots \mathsf{u}_{1}v_{1}$ is in the form (1.1). 
Then for every 
$\bar{\varepsilon} =(\varepsilon _{1},\ldots ,\varepsilon _{n})\in\{ 0, 1\} ^{n}$, we denote by 
$O _{w}^{\bar{\varepsilon}}$ the set 
$\bigcap _{s=1} ^{n}(v_{s}\ldots v_{1}(O_{w}))^{\varepsilon _{s}}$.

\begin{theorem} \label{gab}
Let $G$ act by permutations on some set $X$ and 
$w(\bar{y}) \in \mathbb{F}_{t}\ast G$ be a reduced, non-constant word over $G$ on variables $y_{1},\ldots , y_{t}$. 
Assume that $w(\bar{y}) \in \mathbb{F} _{t}$ or 
$w(\bar{y})=\mathsf{u}_{n}v_{n}\ldots \mathsf{u}_{1}v_{1}$ is in the form (1.1)  
with $O_{w}\neq\emptyset$. 
Assume also that for every  $\bar{\varepsilon}\in\{ 0, 1\} ^{n}$  
\[  
\mbox{ if }O _{w}^{\bar{\varepsilon}}\neq\emptyset\mbox{ , then } 
\mathsf{stab}_{G}(X\setminus O _{w}^{\bar{\varepsilon}})  \mbox{ separates } O _{w}^{\bar{\varepsilon}}.  
\hspace{3cm}(\diamondsuit ) 
\]  
Then the inequality $w(\bar{y}) \neq 1$ has a solution in $G$.
\end{theorem}

\emph{Proof.} 
If $w(\bar{y}) \in\mathbb{F}_{t}$, then $O_{w}=X$ and we simply apply Theorem 1.1 from \cite{A}. 

If $w(\bar{y}) \notin\mathbb{F}_{t}$, then  we follow the proof of Theorem \ref{ab} (keeping the corresponding notation). 
Note that by $(\diamondsuit )$ the set  
$O_{w}^{\bar{\varepsilon}}$ is infinite for 
$O_{w}^{\bar{\varepsilon}}\neq \emptyset$. 
For $O' =O_{w}$ we reformulate Claim M from the proof of Theorem \ref{ab} in the following form.  

\noindent 
 \textbf{Claim M$^{\#}$.} 
 {\em For every $r\leq n$ and every $q\in v_{r}\ldots v_{1}(O_{w})$ there is a unique tuple $\bar{\varepsilon}\in\{ 0, 1\} ^{n}$ such that 
 $q\in O _{w}^{\bar{\varepsilon}}\subseteq v_{r}\ldots v
_{1}(O_{w})$. 
The corresponding $O_{w}^{\bar{\varepsilon}}$ satisfies:} 
\[ 
 \forall V \in \mathcal{V}^{>0}_{w}(O')\, (O _{w}^{\bar{\varepsilon}}\cap V\neq\emptyset\
	\Rightarrow\ O _{w}^{\bar{\varepsilon}}\subseteq V). \hspace{3cm}(\ddag )
	\] 
We prove Claim M$^{\#}$ as follows. 
First observe that 
$\{ O_{w}^{\bar{\varepsilon}}\ |\ \bar{\varepsilon}
\in\{ 0, 1\} ^{n}\}$ 
is a partition of $X$ and hence for any $q\in v_{r}\ldots v_{1}
(O_{w})$ there is a unique tuple
$\bar{\varepsilon}\in\{ 0, 1\} ^{n}$ such that 
$q\in O _{w}^{\bar{\varepsilon}}\subseteq v_{r}\ldots v_{1}
(O_{w})$. 
Now fix this tuple $\bar{\varepsilon}$, take any $s\le n$ and suppose that 
$O _{w}^{\bar{\varepsilon}}\cap v_{s}\ldots v_{1}(O_{w})\neq\emptyset$. 
It follows that $\varepsilon_{s} =1$ in $\bar{\varepsilon}$. 
Thus 
$O _{w}^{\bar{\varepsilon}}\subseteq v_{s}\ldots v_{1}(O_{w})$. 

We now apply the proof of Theorem $\ref{ab}$. 
Let $O'=O_{w}$. 
At the $k$-th step of induction we show that: \parskip2pt
\begin{itemize} 
\item  There is $p\in O'$ and a tuple $\bar{g}=(g_{1},\ldots , g_{t})\in G$ such that $\bar{g}$ is distinctive	for $p$ and $(w)_{k}$. 
\item In the condition above we can choose $\bar{g}$ so that for all $i$, $1\leq i\leq t$,
$\mathsf{supp}(g_{i})\subseteq\mathcal{V}^{>0}_{w}(O')$ and $g_{i}(v_{r}\ldots v_{1}(O'))=v_{r}\ldots v_{1}(O')$ 
for any $r\leq n$. 
\end{itemize} 
At Step 1 we fix any $p\in O_{w}$. 
At Step $k\le L_n$ by Claim M$^{\#}$ we obtain 
$\bar{\varepsilon}\in\{ 0, 1\}^{n}$ such that 
$p_{k-1,\bar{g}} \in O _{w}^{\bar{\varepsilon}}\subseteq v_d \ldots v_{1}(O_{w})$ and $(\ddag )$ holds. 
Hence we may now apply $(\diamondsuit )$ and Remark \ref{rema} to see that for any finite set $Y$ not containing $p_{k-1,\bar{g}}$, the orbit of $p_{k-1,\bar{g}}$ with respect to  
$\mathsf{stab}_{G}((X\setminus O_{w}^{\bar{\varepsilon}})\cup Y)$  is infinite. 

Furthermore, applying the proof of Theorem \ref{ab} we replace each occurance of the the neighborhood 
$O\subseteq v_{d}\ldots v_{1}(O')$ of $p_{k-1,\bar{g}}$ constructed there by the set $O_{w}^{\bar{\varepsilon}}$ as in  Claim M$^{\#}$. 
Consequently, in order to find a desired tuple
$\bar{g}$ we replace the usage of $(\dag )$ from the proof of Theorem \ref{ab} by $(\ddag )$.  
After this modifications the proof of Theorem \ref{ab} works for Theorem \ref{gab}. $\Box$ 

\bigskip 

This theorem cannot be directly applied in the case of subgroups of $S_{\mathsf{fin}}(\mathbb{N})$. 
In this case condition $(\diamondsuit )$ does not hold for $O_w$. 
In Section 3 we will show that this theorem can be applied to so called oscillating identities. 

\subsection{Partial MIF} 

In this section we give a curious application of Theorem \ref{ab}. 
Let $G$ be a group and $E\subseteq G$. 
We denote by $\mathbb{F}_t \ast E$ the set of all reduced words from $(\mathbb{F}_{t}\ast G) \setminus G$ with constants from $E$. 
We remind the reader that $\mathbb{F}_t$ is the free group 
$\mathbb{F}[y_1 ,\ldots , y_t ]$. 
Below we denote by $\mathbb{F}_{\omega}$ the free group 
$\mathbb{F}[z_1 ,\ldots , z_i ,\ldots ]$. 

\begin{definition}
Let $G$ be a group acting on a perfect Hausdorff space $\mathcal{X}$ by homeomorphisms. 
We say that a subset $E\subseteq G$ is  e.o-stable if for every natural $t>0$, each reduced word from $\mathbb{F}_{t}\ast E$ is explicitly oscillating.   
\end{definition} 

\begin{example} \label{MIFTHO} 
{\em Consider Thompson's group $F$ and the subset, say $F_{cf}$, of all $g\in F$ such that $\mathsf{supp}(g)$ is cofinite in $[0,1]$. 
Then $F_{cf}$ is e.o-stable. 
Indeed, when $w(\bar{y})$ is a nontrivial reduced word over $F_{cf}$ then $O_w$ is cofinite in $[0,1]$. 
Note that $F_{cf} \cdot F_{cf} = F$. 
} 
\end{example} 

When a group $G$ is a subgroup of some $\hat{G}$ then it is said that $G$ is {\em existentially closed} in $\hat{G}$ if for every sequence of words 
$$
w_1 (\bar{z}), \ldots , w_{\ell} (\bar{z}), w_{\ell +1} (\bar{z}), \ldots ,w_{\ell +k} (\bar{z}) \in \mathbb{F}_{\omega}\ast G
$$
if the system of equations and inequations 
$$
w_1 (\bar{z})=1, \ldots , w_{\ell} (\bar{z})=1, w_{\ell +1} (\bar{z})\not= 1, \ldots ,w_{\ell +k} (\bar{z})\not= 1
$$ 
has a solution in $\hat{G}$ then it has a solution in $G$. 
It is easy to see that Proposition 5.3 from \cite{HO} can be reformulated as follows.  
\begin{quote} 
A countable group $G$ is MIF if and only if $G$ is existentially closed in  $\mathbb{F}_{t}\ast G$ for every natural $t$. 
\end{quote} 
Moreover, the sufficiency is obvious. 
Thus the following statement is in close relation to it. 

\begin{proposition} \label{e_c} 
Let $G$ be a group acting on a perfect Hausdorff space $\mathcal{X}$ by homeomorphisms and let $t$ be a natural number. 
Assume that a subset $E\subseteq G$ is e.o-stable. 
Then for every sequence of words  
$$
w_1 (\bar{z}), \ldots , w_{\ell} (\bar{z}) , w_{\ell +1} (\bar{z}), \ldots ,w_{\ell +k} (\bar{z}) 
\in \mathbb{F}_{\omega}\ast G 
$$
if the system of equations and inequations 
$$
w_1 (\bar{z})=1, \ldots , w_{\ell} (\bar{z})=1, w_{\ell +1} (\bar{z})\not= 1, \ldots ,w_{\ell +k} (\bar{z})\not= 1
$$ 
has a solution in $\mathbb{F}_t \ast G$ which for every $i\le k$ takes $w_{\ell +i} (\bar{z})$ to some $w'_{\ell +i} (\bar{y}) \in (\mathbb{F}_{t}\ast E)\setminus \{ 1 \}$, then the system has a solution in $G$. 
\end{proposition} 

{\em Proof.}  Let $|\bar{z}| = s$ and let $v_1 (\bar{y}), \ldots , v_s (\bar{y})$ be a solution of the system 
$$
w_1 (\bar{z})=1, \ldots , w_{\ell} (\bar{z})=1, w_{\ell +1} (\bar{z})\not= 1, \ldots ,w_{\ell +k} (\bar{z})\not= 1
$$ 
in $\mathbb{F}_{t}\ast G$ such that  
after the substitution $\bar{z} \leftarrow \bar{v}$ into 
$w_{\ell +1} (\bar{z}), \ldots ,w_{\ell +k} (\bar{z})$ 
we obtain $w'_{\ell +1} (\bar{y}), \ldots ,w'_{\ell +k} (\bar{y})$ from the formulation of the proposition. 
Let 
$$ 
u (y_0,\bar{y}) = [[\ldots [w'_{\ell + 1} (\bar{y}), y^{-1}_0 w'_{\ell +2} (\bar{y}) y_0 ], \ldots ], y^{-k+1}_0 w'_{\ell +k}(\bar{y}) y^{k-1}_0 ] .   
$$ 
We view $u (y_0,\bar{y})$ as an element of $\mathbb{F}[y_0 ,y_{1},\ldots , y_{t}] \ast G$. 
It is easy to see that the normal form of it is non-trivial and belongs to $\mathbb{F}_{t+1}\ast E$, i.e. it is explicitly oscillating.   
Applying Theorem \ref{ab} we find $g_0, g_{1},\ldots , g_{t}\in G$ with $u (g_0,\bar{g})\not=1$. 
Thus 
$w'_{\ell +1} (\bar{g})\not= 1, \ldots ,w'_{\ell +k} (\bar{g})\not=1$. 
In particular the elements $v_1 (\bar{g}), \ldots , v_s (\bar{g})$ form a solution of the system 
$w_{\ell +1} (\bar{z})\not= 1, \ldots ,w_{\ell +k} (\bar{z})\not= 1$ in $G$. 
Since the tuple $v_1 (\bar{y}), \ldots , v_s (\bar{y})$ is a solution of the system 
$w_{1} (\bar{z})= 1, \ldots ,w_{\ell} (\bar{z})= 1$ in $\mathbb{F}_t \ast G$, 
the elements $v_1 (\bar{g}), \ldots , v_s (\bar{g})$ form a solution of the system 
$w_{1} (\bar{z})= 1, \ldots ,w_{\ell} (\bar{z})= 1$ in $G$. 
$\Box$ 

\bigskip 

Using Example \ref{MIFTHO} we see that the statement of this proposition holds when $G$ is Thompson's $F$ and $E=F_{cf}$.

\section{Inequalities in groups with hereditarily separating action}

\subsection{Oscillating words}

In order to apply Theorem \ref{ab} to a broader class of words and a larger number of inequalities we will introduce the notion of an \emph{oscillating} word. 
Intuitively it describes words, which are explicitly oscillating after the transition from the region $O_w$ to some other place. 

Suppose $G$ acts on some Hausdorff topological space $\mathcal{X}$ by homeomorphisms and let $w(\bar{y})$ be a word over $G$ on $t$ variables such that $w(\bar{y})$ is reduced in 
$\mathbb{F} _{t}\ast G$. 
If $w(\bar{y})\notin\mathbb{F} _{t}$ then we write it as follows:  
\[ 
w_{\mathcal{X}}(\bar{y}) =\mathsf{u}_{\mathcal{X},n_{\mathcal{X}}}v_{\mathcal{X},n_{\mathcal{X}}}\mathsf{u}_{\mathcal{X},n_{\mathcal{X}}-1}v_{\mathcal{X},n_{\mathcal{X}}-1}\ldots \mathsf{u}_{\mathcal{X},1}v_{\mathcal{X},1} ,
\] 
where $n_{\mathcal{X}}\in\mathbb{N}$, $\mathsf{u}_{\mathcal{X},i}$ depends only on variables and $v_{\mathcal{X},i}\in G\setminus \{ 1\}$ for $1\le i\leq n_{\mathcal{X}}$. 
Let $v_{\mathcal{X},0}:=1$. 

Suppose that $w_{\mathcal{X}}(\bar{y})$ is not explicitly oscillating. 
We now describe a procedure which produces a family 
$\mathcal{P}^{os}$ of open subsets of $\mathcal{X}$ and a map which associates to each $V\in \mathcal{P}^{os}$ a word  $w_V (\bar{y})$ which is explicitly oscillating in $V$. 
The case $\mathcal{P}^{os}=\emptyset$ is possible but not desirable. 
We call this procedure $\mathsf{Transition}$. 

In the description of it we use the following notation. 
For an open $A\subseteq \mathcal{X}$ let 
$A^{0}:=\mathsf{int}(\mathcal{X}\setminus (A \cup \mathsf{Fix}(G)))$ 
and $A^{1}:=A\setminus \mathsf{Fix}(G)$. 
Below we always assume that $\mathsf{Fix}(G)$ is  finite. 

\bigskip 

\noindent
{\bf Transition.} 
For each sequence 
$\bar{\varepsilon} =(\varepsilon_{1},\ldots , \varepsilon _{n_{\mathcal{X}}})\in\{ 0,1\} ^{n_{\mathcal{X}}}$ we define the set
\[ 
\mathcal{X}_{\bar{\varepsilon}}:=\bigcap _{i=1}^{n_{\mathcal{X}}} v_{\mathcal{X},1}^{-1}\ldots v_{\mathcal{X},i-1}^{-1}\Big( \mathsf{supp}(v_{\mathcal{X},i})
	^{\varepsilon _{i}}\Big)   
\] 
and consider the following family:
\[ 
\mathcal{P}^{1}:=\Big\{ \mathcal{X}_{\bar{\varepsilon}}\ \Big|\ \bar{\varepsilon}\in\{ 0, 1\} ^{n_{\mathcal{X}}}\Big\}\setminus
	\Big\{\emptyset\Big\} . 
\] 
Since $w_{\mathcal{X}} (\bar{y})$ is not explicitly oscillating, 
$\mathcal{X}_{(1,\ldots , 1)}=\emptyset$, i.e. not in $\mathcal{P}^{1}$. 
For each $\mathcal{X}_{\bar{\varepsilon}}\in\mathcal{P}^{1}$ we define a word $w'_{\mathcal{X}_{\bar{\varepsilon}}}(\bar{y})$ in the following way:
\[ 
w'_{\mathcal{X}_{\bar{\varepsilon}}}(\bar{y}):=\mathsf{u}_{\mathcal{X},n_{\mathcal{X}}}v^{\varepsilon _{n_{\mathcal{X}}}}_{\mathcal{X},n_{\mathcal{X}}}\mathsf{u}_{\mathcal{X},n_{\mathcal{X}}-1}v^{\varepsilon _{n_{\mathcal{X}}-1}}_{\mathcal{X},n_{\mathcal{X}}-1}\ldots \mathsf{u}_{\mathcal{X},1}v^{\varepsilon _{1}}_{\mathcal{X},1}.  
\] 
Then we reduce $w'_{\mathcal{X}_{\bar{\varepsilon}}}(\bar{y})$ in $\mathbb{F}_{t}\ast G$. 
If after the reduction the obtained word $w'(\bar{y})$ is of the form 
$\, \, \cdot \cdot \cdot$ $v'\mathsf{u}'$,
where $\mathsf{u}'$ contains only variables and $v'\in G\setminus\{ 1\}$  then we conjugate $w'(\bar{y})$ by $(\mathsf{u}')^{-1}$ and denote the obtained word by $w_{\mathcal{X}_{\bar{\varepsilon}}}(\bar{y})$. 
Otherwise we simply take $w_{\mathcal{X}_{\bar{\varepsilon}}}(\bar{y}):=w'(\bar{y})$. 
Let 
\[ 
\mathcal{W}^{1} = \Big\{ w_{\mathcal{X}_{\bar{\varepsilon}}}(\bar{y}) \ \Big| \   w_{\mathcal{X}_{\bar{\varepsilon}}}(\bar{y})\not= 1 \mbox{ , }  
\mathcal{X}_{\bar{\varepsilon}}\in\mathcal{P} ^{1} \Big\}. 
\]  
If there is some word 
$w_{\mathcal{X}_{\bar{\varepsilon}}}(\bar{y})\in\mathcal{W} ^{1}$ which is explicitly oscillating in $\mathcal{X}_{\bar{\varepsilon}}$, let
\[ 
\mathcal{P}^{os}:=\Big\{ \mathcal{X}_{\bar{\varepsilon}}\in\mathcal{P} ^{1}\ \Big|\ w_{\mathcal{X}_{\varepsilon}}(\bar{y}) 
	\mbox{ is explicitly oscillating in }\mathcal{X}_{\bar{\varepsilon}}\Big\} . 
\] 
If there is no explicitly oscillating 
$w_{\mathcal{X}_{\bar{\varepsilon}}}(\bar{y})\in\mathcal{W} ^{1}$ then for every $\mathcal{X}_{\bar{\varepsilon}}\in\mathcal{P} ^{1}$ we repeat the process described above replacing $\mathcal{X}$ by $\mathcal{X}_{\bar{\varepsilon}}$ and $w_{\mathcal{X}}(\bar{y})$ by $w_{\mathcal{X}_{\bar{\varepsilon}}}(\bar{y})$. 
For every $\mathcal{X}_{\bar{\varepsilon}}\in\mathcal{P} ^{1}$ we define the family $\mathcal{P} ^{2} _{\mathcal{X}_{\bar{\varepsilon}}}$ and the
corresponding set of words 
$\mathcal{W} ^{2} _{\mathcal{X}_{\bar{\varepsilon}}}$ exactly as 
$\mathcal{P} ^{1}$ and $\mathcal{W} ^{1}$ were defined above. 
Now let
\[ 
\mathcal{P} ^{2}:=\bigcup\Big\{ \mathcal{P} ^{2} _{\mathcal{X}_{\bar{\varepsilon}}}\ \Big|\ \mathcal{X}_{\bar{\varepsilon}}\in \mathcal{P} ^{1}\Big\} . 
\] 
Let $\mathcal{W} ^{2}$ be the set of all words $w_{V}(\bar{y})$ for $V\in\mathcal{P}^{2}$ defined as 
$w_{\mathcal{X}_{\bar{\varepsilon}}}(\bar{y})$ above. 
If there is an explicitly oscillating word in $\mathcal{W}^{2}$ then we define
\[ 
\mathcal{P} ^{os}:=\Big\{ V\in\mathcal{P} ^{2}\ \Big|\ w_{V}(\bar{y}) \mbox{ is explicitly oscillating in }V\Big\} 
\] 
and finish the construction. 
If there is no explicitly oscillating word in $\mathcal{W} ^{2}$, then we continue this procedure. 
If for some $k\in\mathbb{N}$, $\mathcal{W}^{k}$ contains explicitly oscillating words or $\mathcal{W}^{k}=\emptyset$, then the procedure terminates. 

To unify our notation we will denote $\mathcal{W}^{0}:=\{ w_{\mathcal{X}}(\bar{y})\}$. 

\begin{lemma}
Procedure $\mathsf{Transition}$ terminates after finitely many steps.
\end{lemma}

\emph{Proof.} 
Suppose that for some $k>1$ we are given a word 
$w_{V}(\bar{y})\in\mathcal{W}^{k-1}$,
which is not explicitly oscillating. 
Thus $w_V (\bar{y})\not\in\mathbb{F} _{t}$, i.e. $w_V (\bar{y})$  contains some constants $v_{V,i}\in G$, $1\leq i\leq n_{V}$. 

Let 
$\bar{\varepsilon}\in\{ 0, 1\}^{n_{V}}$, $V_{\bar{\varepsilon}}\in\mathcal{P}^{k} _{V}$ and 
$w_{V_{\bar{\varepsilon}}}(\bar{y})\in\mathcal{W}^{k}$ be obtained from $w_{V}(\bar{y})$ as in the construction.
Since $w_{V}(\bar{y})$ was not explicitly oscillating, 
\[ 
U_{1,\ldots ,1}=\bigcap _{i=1}^{n_{V}} v_{V,1}^{-1}\ldots v_{V,i-1}^{-1}\Big( \mathsf{supp}(v_{V,i})^{1} \Big) =\emptyset. 
\]
Thus for some $i$, $1\leq i\leq n_{V}$, $v_{V, i}$ becomes $\mathsf{id}$ in the word $w_{V_{\bar{\varepsilon}}}(\bar{y})$. 
Therefore the number of constants in $w_{V_{\bar{\varepsilon}}}(\bar{y})$ is strictly smaller than the corresponding number in  $w_{V}(\bar{y})$. 

We now see that after finitely many steps we find some 
$k$ such that either $\mathcal{W}^{k}$ contains an explicitly oscillating word or the words 
$w_{V}(\bar{y})$ are equal to $1$ for all $V\in\mathcal{P}^{k}$.  $\Box$ 

\bigskip 

\begin{definition} \label{trans}
Under the notation above the initial word $w_{\mathcal{X}}(\bar{y})$ is called oscillating if it is explicitly oscillating or procedure $\mathsf{Transition}$ terminates producing the set $\mathcal{P} ^{os} \not=\emptyset$.  
If the initial word $w_{\mathcal{X}}(\bar{y})$ is not oscillating and  procedure $\mathsf{Transition}$ teminates for some $k\in\mathbb{N}$ such that for all  $V\in\mathcal{P} ^{k}$, $w_{V}(\bar{y})=1$, then the word $w_{\mathcal{X}}(\bar{y})$ is called rigid. 
\end{definition} 

Assume that $w_{\mathcal{X}}(\bar{y})$ is oscillating but is not explicitly oscillating.  
If $V\in\mathcal{P}^{os}$ and $w_{V}\notin\mathbb{F}_{t}$ (considered in $\mathbb{F}_{t}\ast G$), we view it in the form 
$w_{V}(\bar{y})=\mathsf{u}_{V,n_{V}}v_{V,n_{V}}\ldots \mathsf{u}_{V,1}v_{V,1}$, where 
$n_{V}\in\mathbb{N}$ and for each $i\leq n_{V}$ the subword $\mathsf{u}_{V,i}$ depends only on variables and 
$v_{V,i}\in G\setminus\{ 1\}$. 
Let $v_{V,0}=1$. 

Define
\[ 
\hat{O}_{w}:=\bigcup _{V\in\mathcal{P} ^{os}} \Big( (V\setminus \mathsf{Fix}(G))\cap\bigcap _{i=0} ^{n_{V}-1} v_{V,0}^{-1}v
_{V,1}^{-1}\ldots v_{V,i}^{-1}\Big( \mathsf{supp}(v_{V,i+1})\Big)\Big) .  
\] 
In particular, for  $w_{V}(\bar{y})\in \mathbb{F}_{t}$, we have $n_{V}=0$ and the contribution of $w_{V}(\bar{y})$ to the above union equals $V\setminus \mathsf{Fix}(G)$. 

To unify notation in the case of explicitly oscillating $w_{\mathcal{X}} (\bar{y})$ we put $\hat{O}_w = O_w$.  

\bigskip

\begin{example} \label{ex1}
{\em Consider Thompson's group $F$ with its standard action on $[0,1]$. 
We now give several illustrations of notions introduced above. 
There are graphic illustration of these cases in \cite{Zarphd}. 
} 

\noindent 
(a) 
{\em Let 

	$\bullet$ \ \ $w_{2}(y)=x_{[0,\frac{1}{2}],0}^{-1}yx_{[\frac{1}{2},1],1}^{-1}y^{-1}x_{[0,\frac{1}{2}],1} yx_{[0,\frac{1}{2}],2}^{-1}$. \\ 
The word $w_{2}(y)$ is not explicitly oscillating.  
Indeed, 
$x_{[0,\frac{1}{2}],2}x_{[0,\frac{1}{2}],1}^{-1}((\frac{1}{2},1))\cap (0,\frac{1}{2})=\emptyset$.  
We state that it is oscillating. 
To see this we apply $\mathsf{Transition}$.   
There are four constant segments in $w_{2}(y)$: 
$v_{[0,1],1}=x_{[0,\frac{1}{2}],2}^{-1}$ ,
$v_{[0,1],2}=x_{[0,\frac{1}{2}],1}$ , 
$v_{[0,1],3}=x_{[\frac{1}{2},1],1}^{-1}$ and 
$v_{[0,1],4}=x_{[0,\frac{1}{2}],0}^{-1}$. 
It is clear that \\

$\ \ \ \mathsf{supp}(v_{[0,1],1})^{1}=\Big(\frac{3}{8},\frac{1}{2}\Big)$ and 
$\mathsf{supp}(v_{[0,1],1})^{0}=\Big( 0,\frac{3}{8}\Big)\cup\Big(\frac{1}{2},1\Big)$,\\

$\ \ \ \mathsf{supp}(v_{[0,1],2})^{1}=\Big(\frac{1}{4},\frac{1}{2}\Big)$ and  
$\mathsf{supp}(v_{[0,1],2})^{0}= \Big( 0,\frac{1}{4}\Big)\cup\Big(\frac{1}{2},1\Big)$,\\

	$\ \ \ \mathsf{supp}(v_{[0,1],3})^{1}=\Big(\frac{3}{4},1\Big)$ and  
$\mathsf{supp}(v_{[0,1],3})^{0}=\Big( 0, \frac{3}{4}\Big)$,\\

	$\ \ \ \mathsf{supp}(v_{[0,1],4})^{1}=\Big( 0,\frac{1}{2}\Big)$ and 
$\mathsf{supp}(v_{[0,1],4})^{0}=\Big( \frac{1}{2},1\Big)$.\\

\noindent 
Thus the family 
$\mathcal{P}^{1}$ for $w_{2}(y)$ equals 
$\{ (0,\frac{1}{4}), (\frac{1}{4},\frac{3}{8}), (\frac{3}{8},\frac{1}{2}), (\frac{1}{2},\frac{3}{4}),(\frac{3}{4},1)\}$. 
Hence we obtain five reduced words: \\ 
$(w_2)'_{(0,\frac{1}{4})}(y)=x_{[0,\frac{1}{2}],0}^{-1}y$ , 
\hspace{2cm} $(w_2)'_{(\frac{1}{4}, \frac{3}{8})}(y)=x_{[0,\frac{1}{2}],0}^{-1}x_{[0,\frac{1}{2}],1}y$ , \\  
$(w_2)'_{(\frac{3}{8},\frac{1}{2})}(y) =x_{[0,\frac{1}{2}],0}^{-1}x_{[0,\frac{1}{2}],1}yx_{[0,\frac{1}{2}],2}^{-1}$ , $(w_2)'_{(\frac{1}{2},\frac{3}{4})}(y)=y$ 
\hspace{1cm} and \\ 
$(w_2)'_{(\frac{3}{4},1)}(y)=yx_{[\frac{1}{2},1],1}^{-1}$. \\ 
The corresponding words  
$(w_2)_{(0,\frac{1}{4})}(y)$, $(w_2)_{(\frac{1}{4},\frac{3}{8})}(y)$, $(w_2)_{(\frac{1}{2},\frac{3}{4})}(y)$ and $(w_2)_{(\frac{3}{4},1)}(y)$ are explicitly oscillating. 
Note that $(w_2)_{(\frac{3}{8},\frac{1}{2})}(y)$ is non-trivial and not explicitly oscillating in 
$(\frac{3}{8},\frac{1}{2})$.  
Indeed, 
$\mathsf{supp}(x_{[0,\frac{1}{2}],0}^{-1}x_{[0,\frac{1}{2}],1})\cap \mathsf{supp}(x_{[0,\frac{1}{2}],2}^{-1})=\emptyset$. 
This finishes the procedure and we see that
$w_{2}(y)$ is oscillating, where  
$\hat{O}_{w_2} = (0,\frac{3}{8})\cup (\frac{1}{2},1)$. 
} 

\noindent 
(b) 
{\em The following word is explicitly oscillating and has trivial product of constants: 

	$\bullet$ \ \ $w_3(y) =yx_{1}y^{-1}x_{1}^{-1}.$ 

\noindent 	
Indeed, 
\[ 
O_{w_{3}}=x_{1}\Big(\Big(\frac{1}{2},1\Big)\Big)\cap\Big(\frac{1}{2},1\Big) =\Big(\frac{1}{2},1\Big) ,  
\] 
and $x_{1} x_{1}^{-1}=\mathsf{id}$. \\ 
On the other hand the word 

	$\bullet$ \ \ $w_4 (y) =yx_{1}y^{-1}x_{[0,\frac{1}{2}],0}y^{2}x_{1}^{-1}$ \\ 
is not explicitly oscillating, because 
\[ 
O_{w_{4}}(y) =x_{1}x_{[0,\frac{1}{2}],0}^{-1}\Big(\Big(\frac{1}{2},1\Big)\Big)\cap
x_{1}\Big(\Big( 0,\frac{1}{2}\Big)\Big)\cap\Big(\frac{1}{2},1\Big) = \emptyset.
\] 
Since $x_{1}x_{[0,\frac{1}{2}],0}x_{1}^{-1}=x_{[0,\frac{1}{2}],0}$ and
$\mathsf{supp}(x_{[0,\frac{1}{2}],0})=(0,\frac{1}{2})$, 
the word $w_4 (y)$ has non-trivial product of constants. 
We will see in Proposition \ref{con-tra} that this guarantees the $w_4 (y)$ is oscillating. 
The corresponding pictures can be found in \cite{Zarphd}.  
} 
\end{example} 

\bigskip 

\begin{example} 
{\em We start with a word, denoted by $w_5 (y)$, which is rigid under the standard action of Thompson's group $F$.  
Using this we construct $w_6 (y_1 , y_2 )$ such that $\mathsf{Transition}$ needs 2 steps in order to show that it is oscillating  } 

\noindent 
(a) 
{\em Let 

	$\bullet$ \ \ $w_{5}(y)=y^{-1}x_{1}yx_{[0,\frac{1}{2}],0}y^{-1}x_{1}^{-1}yx_{[0,\frac{1}{2}],0}^{-1}$. \\ 
This word has four constant segments: 
$v_{[0,1],1}=x_{[0,\frac{1}{2}],0}^{-1}$ , 
$v_{[0,1],2}=x_{1}^{-1}$ , 
$v_{[0,1],3}=x_{[0,\frac{1}{2}],0}$ , 
$v_{[0,1],4}=x_{1}$ , which define the following supports: \\

	$\ \ \ \mathsf{supp}(v_{[0,1],1})^{1}=\Big( 0,\frac{1}{2}\Big)$  and 
$\mathsf{supp}(v_{[0,1],1})^{0}=\Big(\frac{1}{2},1\Big)$,\\

	$\ \ \ \mathsf{supp}(v_{[0,1],2})^{1}=\Big(\frac{1}{2},1\Big)$  and  
$\mathsf{supp}(v_{[0,1],2})^{0}=\Big( 0, \frac{1}{2}\Big)$.\\

	$\ \ \ \mathsf{supp}(v_{[0,1],3})^{1}=\Big( 0,\frac{1}{2}\Big)$  and  
$\mathsf{supp}(v_{[0,1],3})^{0}=\Big(\frac{1}{2}, 1\Big)$,\\

	$\ \ \ \mathsf{supp}(v_{[0,1],4})^{1}=\Big(\frac{1}{2},1\Big)$ and  
$\mathsf{supp}(v_{[0,1],4})^{0}=\Big(0 ,\frac{1}{2}\Big)$.\\

\noindent 
The family $\mathcal{P}^{1}$ for $w_5(y)$ equals 
$\{ (0,\frac{1}{2}),(\frac{1}{2},1)\}$.
	Thus $\mathsf{Transition}$ gives two words: 
\[ 
(w_5)'_{(0,\frac{1}{2})}(y)=y^{-1}yx_{[0,\frac{1}{2}],0}y^{-1}yx_{[0,\frac{1}{2}],0}^{-1} 
\] 
and 
\[ 
(w_5)'_{(\frac{1}{2},1)}(y)=y^{-1}x_{1}yy^{-1}x_{1}^{-1}y . 
\] 
In fact both of them are equal to $\mathsf{id}$. 
Hence $\mathcal{W}^{1}$ is empty and therefore $w_5 (y)$ is rigid. }
 
\noindent 
(b) {\em In this example we show that $\mathsf{Transition}$ sometimes needs several steps. 
Let $v$ and $v'$ be elements of $F$ with $\mathsf{supp}(v) = [0,1]= \mathsf{supp}(v')$. 
When $[a,b]$ is a dyadic subinterval of $[0,1]$ let $v_{[a,b]}$ and $v'_{[a,b]}$ be the corresponding elements defined on $[a,b]$.  
Let \\ 
$v_{[0,1],1} = v'_{[0,\frac{1}{4}]}x_{[\frac{1}{4},\frac{1}{2}],0}x_{[\frac{1}{2},\frac{3}{4}],0}v_{[\frac{3}{4}, 1]}$ \, 
\, , \, \, $v_{[0,1],6} = v'_{[0,\frac{1}{4}]}x^{-1}_{[\frac{1}{4},\frac{1}{2}],0}x^{-1}_{[\frac{1}{2},\frac{3}{4}],0}v_{[\frac{3}{4}, 1]}$ , \\
$v_{[0,1],7} = x_{[0,\frac{1}{4}],0}v'_{[\frac{1}{4},\frac{1}{2}]}v_{[\frac{1}{2},\frac{3}{4}]}x_{[\frac{3}{4},1],0}$
\, and \, 
$v_{[0,1],12} = x^{-1}_{[0,\frac{1}{4}],0}v'_{[\frac{1}{4},\frac{1}{2}]}v_{[\frac{1}{2},\frac{3}{4}]}x^{-1}_{[\frac{3}{4},1],0}$. \\ 
Note that these elements of $F$ have the same support: 
$(0,1) \setminus \{ \frac{1}{4}, \frac{1}{2}, \frac{3}{4} \}$. 
Thus computing $\mathsf{supp}^0$ of them we obtain $\emptyset$ in each case. 
In order to define the remaining constants we use $\frac{1}{2}$-versions of the word $w_5$ defined in (a).  
Let \\ 
$v_{[0,1],2} = v_{[0,1],8} = x^{-1}_{[0,\frac{1}{4}],0}x^{-1}_{[\frac{1}{2},\frac{3}{4}],0}$ \, \, , \, \,  
$v_{[0,1],3} = v_{[0,1],9} = x^{-1}_{[0,\frac{1}{2}],1}x^{-1}_{[\frac{1}{2},1],1}$ , \\  
$v_{[0,1],4} = v_{[0,1],10} = x_{[0,\frac{1}{4}],0}x_{[\frac{1}{2},\frac{3}{4}],0}$ \, and \, 
$v_{[0,1],5} = v_{[0,1],11} = x_{[0,\frac{1}{2}],1}x_{[\frac{1}{2},1],1}$ . \\ 
Let 
\[ 
w_6 = v_{12}y^{-1}_2 y^{-1}_1 v_{11} y_1 v_{10} y^{-1}_1 v_9 y_1 v_8 y_2 v_7 y_1 y_2 v_6 y^{-1}_2 y^{-1}_1 v_{5} y_1 v_{4} y^{-1}_1 v_3 y_1 v_2 y_2 v_1 . 
\]  
We omit $[0,1]$-indexes. 
Note that between $v_{12}y^{-1}_2$ and $y_2 v_7$ we have $[0, \frac{1}{2}]$- and $[\frac{1}{2},1]$-versions of $w_5$ placed in parallel. 
The same word appears between $v_{6}y^{-1}_2$ and $y_2 v_1$. 
The family of open subsets of $[0,1]$ computed by $\mathsf{Transition}$ at step 1 is as follows: $\mathcal{P}^1 = \{ (0, \frac{1}{4}) \cup (\frac{1}{2}, \frac{3}{4} )$ , 
$ (\frac{1}{4} , \frac{1}{2} )\cup (\frac{3}{4}, 1 ) \}$.  
Applying the analysis of case (a) we obtain that the subword of $w_6$ between $v_{12}y^{-1}_2$ and $y_2 v_7$ is is trivial for each element of $\mathcal{P}^1$. 
The same argument shows that the subword of $w_6$ between $v_{6}y^{-1}_2$ and $y_2 v_1$ is trivial too.  
After cancellation of appearing subwords $y^{-1}_2 y_2$ we have that $w_6$ becomes $v_{12}v_7 y_1 y_2 v_6 v_1$ on each element of $\mathcal{P}^1$. 
Note that $\mathsf{supp}(v_{12}v_7 ) = ( \frac{1}{4} ,\frac{3}{4}) \setminus \{ \frac{1}{2} \}$ and 
$\mathsf{supp}(v_{6}v_1 ) = (0, \frac{1}{4} )\cup (\frac{3}{4}, 1 )$. 
In particular the word $v_{12}v_7 y_1 y_2 v_6 v_1$ is not explicitly oscillating and the open sets mentioned in the previous sentence form $\mathcal{P}^2$. 
It is easy to see that 
$(w_6)_{\mathsf{supp}(v_{12}v_7 )} = y^{-1}_1y^{-1}_2 v_{12}v_7$ and 
$(w_6)_{\mathsf{supp}(v_6v_1 )} = y_1y_2 v_6v_1$. 
These words are explicitly oscillating. 
}
\end{example} 

\bigskip 

The following lemma exhibits correspondence between the existence of solutions of the inequality $w_{\mathcal{X}} (\bar{y})\neq 1$ and $w_{V}(\bar{y})\neq 1$, where $w_{V}(\bar{y})$ is derived from $w_{\mathcal{X}}(\bar{y})$ by $\mathsf{Transition}$. 
We use the notation of this section and $\mathsf{Transition}$. 

\begin{lemma}\label{ii}
Suppose that $k\geq 1$, $U\in\mathcal{P}^{k}$, $p\in U$ 
and $\bar{g}=(g_{1},\ldots , g_{t})\in G$. 
Assume $w_{U}(\bar{g})(p)\neq p$ where 
$w_{U}(\bar{y}) \in\mathcal{W}^{k}$ corresponds to $U$. 
If for every $i$, $1\leq i\leq t$, and every 
$(\varepsilon_1 ,\ldots ,\varepsilon_j )\in \{ 0,1\}^j$,
$1\leq j\leq n_{\mathcal{X}}$, the element $g_{i}$ stabilizes $v^{\varepsilon_1}_{\mathcal{X},j}\ldots v^{\varepsilon_j}_{\mathcal{X},1}(U)$ setwise, then $w_\mathcal{X} (\bar{g})\neq 1$.
\end{lemma}

\emph{Proof.} 
Fix some $\bar{g}\in G$ satisfying the conditions of the lemma. 
Let $r\leq k$, $U\subseteq V\subseteq V'$, 
$V\in\mathcal{P}^{r}_{V'}$ and $w_{V'}(\bar{y})$ be a word from $\mathcal{W}^{r-1}$ which is not explicitly oscillating. 
Consider the word $w_{V}(\bar{y})\in\mathcal{W}^{r}$, which 
is obtained from $w_{V'}(\bar{y})$ by the appropriate reductions and conjugation. 
We will show that, if for some $p'\in U$, $w_{V}(\bar{g})(p')\neq p'$, then we have $w_{V'}(\bar{g})(p'')\neq p''$, where $p''\in U$ is obtained from $p'$ by (possibly) several applications of $g_i$, $1\le i \le t$.
This proves the lemma by induction starting with the case $r=k$ and $w_{U}(\bar{y})\in\mathcal{W}^{k}$, 
where 
$p' := p\in U$ as in the formulation. 
For ease of notation we will put $p'=p$ below. 
Let 
\[ 
w_{V'}(\bar{y}):=\mathsf{u}_{V',n_{V'}}v_{V',n_{V'}}\mathsf{u}_{V',n_{V'}-1}v_{V',n_{V'}-1}\ldots \mathsf{u}_{V',1}v_{V',1}, 
\] 
where $n_{V'}\in\mathbb{N}\setminus\{ 0\}$, $\mathsf{u}_{V',i}$ depends only on variables and $v_{V',i}\in G\setminus\{1\}$ for $i\leq n_{V'}$. 
Assume that $w_{V}(\bar{y})$ is obtained by reductions and conjugation from the word
\[ 
w'_{V'_{\bar{\varepsilon}}}(\bar{y})=\mathsf{u}_{V',n_{V'}}v^{\varepsilon _{n_{V'}}}_{V',n_{V'}}\mathsf{u}_{V',n_{V'}-1} v^{\varepsilon _{n_{V'}-1}}_{V',n_{V'}-1}\ldots \mathsf{u}_{V',1}v^{\varepsilon _{1}}_{V',1}, 
\]
for some $\bar{\varepsilon}\in\{ 0, 1\}^{n_{V'}}$. 
Note that from the description of $\mathsf{Transition}$ and the assumption that every element $g_{i}$ stabilizes every $v^{\varepsilon_1}_{\mathcal{X},j}\ldots v^{\varepsilon_j}_{\mathcal{X},1}(U)$ setwise we have that for every $j\leq n_{V'}$, elements of $\bar{g}$ stabilize $v_{V',j}\ldots v_{V',1}(U)$ setwise. 
At this step of induction we additionally assume 
that in order to produce $w_V(\bar{y})$ from  $w_{V'}(\bar{y})$ conjugation is not used. 
Otherwise we replace $p$ by an appropriate $\mathsf{u}'(\bar{g})^{-1}(p)$ which is still in $U$. 

Let us compute $w_{V'}(\bar{g})(p)$. 
To simplify notation for each
$j\leq n_{V'}$ denote by $p'_{j}$ the point $\mathsf{u}_{V',j-1}(\bar{g})v_{V',j-1}\ldots \mathsf{u}_{V',1}(\bar{g})
v_{V',1}(p)$ and by $p_{j}$ the point 
\[ 
\mathsf{u}_{V',j-1}(\bar{g})v^{\varepsilon _{j-1}}_{V',j-1}\ldots \mathsf{u}_{V',1}(\bar{g})v^{\varepsilon_{1}}_{V',1}(p). 
\] 
We put $p'_{1}=p_{1} =p$ and 
prove by induction that $p_{j}=p_{j}'$. 

\medskip 
\noindent
 \textbf{Claim.} 
Assume $j\leq n_{V'}$, and for each $i\le j$ the points $p_i$ and $p'_i$ coincide.  
Then 
\[ 
\mathsf{u}_{V',j}(\bar{g})v_{V',j}(p'_{j})=\mathsf{u}_{V',j}(\bar{g})v^{\varepsilon _{j}}_{V',j}(p_{j}). \hspace{3cm}(\dag )
\] 
Note that in the case $j=n_{V'}$ the claim implies the following in/equality 
\[ 
w_{V'}(\bar{g})(p)=\mathsf{u}_{V',n_{V'}}(\bar{g})v_{V',n_{V'}}(p_{n_{V'}}')=w_{V}(\bar{g})(p)\neq p. 
\] 
This will finish the proof of the lemma. 

\medskip 
\noindent 
\emph{Proof of the claim.} 
If $\varepsilon _{j}=1$, then 
$\mathsf{u}_{V',j}(\bar{g})v_{V',j}=\mathsf{u}_{V',j}(\bar{g})v^{\varepsilon _{j}}_{V',j}$ and we are done. 

 Now assume that $\varepsilon _{j}=0$. 
We claim that $p'_{j}\in \mathcal{X}\setminus \mathsf{supp}(v_{V',j})$. 
Indeed, since
\[ 
V=V'_{\bar{\varepsilon}}=\bigcap _{i=1}^{n_{V'}} v_{V',1}^{-1}\ldots v_{V',i-1}^{-1}
\Big( \mathsf{supp}(v_{V',i})^{\varepsilon _{i}}\Big) , 
\] 
for each $\ell\leq n_{V'}-1$ either 
$v_{V',\ell}\ldots v_{V',1}(V)\subseteq \mathsf{supp}(v_{V',\ell +1})$ 
or 
$v_{V',\ell}\ldots v_{V',1}(V)\subseteq \mathcal{X}\setminus \mathsf{supp}(v_{V',\ell +1})$. 
Since for every $i\leq t$ and every $\ell \leq n_{V'}$ the element $g_{i}$ stabilizes
$v_{V',\ell}\ldots v_{V',1}(U)$, we see 
\[ 
\mathsf{u}_{V',j-1}(\bar{g})v_{V',j-1}\ldots \mathsf{u}_{V',1}(\bar{g})v_{V',1}(U)\subseteq\Big( \mathcal{X}\setminus \mathsf{supp}(v_{V',j})\Big) . 
\] 
Since $p\in U$, 
we have $p'_{j}\in \mathcal{X}\setminus \mathsf{supp}(v_{V',j})$ and 
\[ 
u_{V',j}(\bar{g})v_{V',j}(p'_{j})=u_{V',j}(\bar{g})(p'_{j})=u_{V',j}(\bar{g})v^{\varepsilon _{j}}
_{V',j}(p_{j}). 
\] 
$\square$ 

\bigskip 

\begin{example} 
{\em 
Consider any $G \le S_{\mathsf{fin}}(\mathbb{N})$ with respect to the action on $\mathcal{X} =\mathbb{N}$. 
Assume that the word 
\[ 
w(\bar{y}) =\mathsf{u}_{n}v_{n}\mathsf{u}_{k-1}v_{k-1}\ldots \mathsf{u}_{1}v_{1} ,
\] 
is given in the form (1.1) and the corresponding word 
$\mathsf{u}_{n}\mathsf{u}_{n-1}\ldots \mathsf{u}_{1}$ 
(after reductions of all $v_i$) is non-trivial.  
Then $w(\bar{y})$ is oscillating. 
Indeed, applying the first step of $\mathsf{Transition}$ we see that the region  
\[ 
\mathcal{X}_{\bar{0}}:=\bigcap _{i=1}^{n_{\mathcal{X}}} v_{\mathcal{X},1}^{-1}\ldots v_{\mathcal{X},i-1}^{-1}\Big( \mathsf{supp}(v_{\mathcal{X},i})^{0}\Big)  
\] 
belongs to $\mathcal{P}^{os}$. 
As we already know $S_{\mathsf{fin}}(\mathbb{N})$ is hereditarily separating  on $\mathbb{N}$. 
By Abert's theorem from \cite{A} the inequality $\mathsf{u}_{k}\mathsf{u}_{k-1}\ldots \mathsf{u}_{1}\not=1$ has a solution in $S_{\mathsf{fin}}(\mathbb{N})$. 
By Lemma \ref{ii} we obtain a solution of $w(\bar{y}) \not=1$. 
This observation gives some kind of extension of Theorem \ref{gab}.}
\end{example} 

\bigskip 

The following proposition shows that oscillation unifies explicit oscillation with non-triviality of product of constants. 

\begin{proposition} \label{con-tra}
Assume that $G$ acts on a Hausdorff topological space $\mathcal{X}$ by homeomorphisms. 
Let $w(\bar{y}) \in \mathbb{F}_t * G$ be a word in the form (1.1). 
If $w(\bar{y})$ has non-trivial product of constants, then $w(\bar{y})$ is oscillating.  
\end{proposition} 

{\em Proof.} 
Assume that $w(\bar{y})$ is not explicitly oscillating and is in the form as in the beginning of the section: 
\[ 
w_{\mathcal{X}}(\bar{y}) =\mathsf{u}_{\mathcal{X},n_{\mathcal{X}}}v_{\mathcal{X},n_{\mathcal{X}}}\mathsf{u}_{\mathcal{X},n_{\mathcal{X}}-1}v_{\mathcal{X},n_{\mathcal{X}}-1}\ldots \mathsf{u}_{\mathcal{X},1}v_{\mathcal{X},1} ,
\] 
and $\mathsf{supp}(v_{\mathcal{X},n_{\mathcal{X}}}v_{\mathcal{X},n_{\mathcal{X}}-1}\ldots v_{\mathcal{X},1})$ 
is a non-empty subset of $\mathcal{X}$.  
Let $p$ belong to this support. 
Let $p_j = v_{\mathcal{X},j}\ldots v_{\mathcal{X},1}(p)$, $1 \le j \le n_\mathcal{X}$.  
We may assume that when the point $p_j$ belongs to $\overline{\mathsf{supp}(v_{\mathcal{X},j})}$ then it already belongs to $\mathsf{supp}(v_{\mathcal{X},j})$. 
Indeed, if $p_j\in \overline{\mathsf{supp}(v_{\mathcal{X},j})}\setminus\mathsf{supp}(v_{\mathcal{X},j})$ 
then we replace $p$ by a sufficiently close $p'\in \mathsf{supp}(v_{\mathcal{X},n_{\mathcal{X}}}v_{\mathcal{X},n_{\mathcal{X}}-1}\ldots v_{\mathcal{X},1})$ so that the corresponding $p'_j$ belongs to $\mathsf{supp}(v_{\mathcal{X},j})$. 
Using continuity of all $v_{\mathcal{X},i}$ 
we arrange that for each $i\not= j$  
\[ 
p_i \in \mathsf{supp}(v_{\mathcal{X},i}) \Rightarrow 
p'_i \in \mathsf{supp}(v_{\mathcal{X},i}) 
\mbox{ \, and \, } 
p_i \not\in \overline{\mathsf{supp}(v_{\mathcal{X},i})} \Rightarrow 
p'_i \not\in \overline{\mathsf{supp}(v_{\mathcal{X},i})}. 
\] 
Repeating such replacement several times we eventually obtain the required property. 

For each $j\le n_{\mathcal{X}}$ we define $\varepsilon_j$ to be $1$ if  $p_j \not= p_{j-1}$ and to be $0$ otherwise.   
In particular we see that when $\varepsilon_j = 0$, the point $p_j$ does not belong to $\overline{\mathsf{supp}(v_{\mathcal{X},j})}$.   
Then it is easy to see that  
$p\in \mathcal{X}_{\bar{\varepsilon}}$ where 
\[ 
\mathcal{X}_{\bar{\varepsilon}}:=\bigcap _{i=1}^{n_{\mathcal{X}}} v_{\mathcal{X},1}^{-1}\ldots v_{\mathcal{X},i-1}^{-1}\Big( \mathsf{supp}(v_{\mathcal{X},i})^{\varepsilon _{i}}\Big) . 
\] 
Let $w_{\mathcal{X}_{\bar{\varepsilon}}}(\bar{y})$ be the corresponding word obtained by $\mathsf{Transition}$. 
It has fewer constants than $w(\bar{y})$ has, by the assumption that the latter is not explicitly oscillating. 
Furthermore, 
$w_{\mathcal{X}_{\bar{\varepsilon}}}(\bar{1})(p) = w_{\mathcal{X}}(\bar{1})(p)$ and 
$p$ obviously belongs to the support of the product of constants of $w_{\mathcal{X}_{\bar{\varepsilon}}}(\bar{y})$. 

We now apply the same procedure to 
$w_{\mathcal{X}_{\bar{\varepsilon}}}(\bar{y})$ and iterate it until we obtain an explicitly oscillating word. 
The necessity of the latter output follows from the fact that 
at each step we obtain a shorter word with non-trivial product of constants.   
$\Box$ 

\subsection{Solving a system of inequalities in the case of oscillating words}

In this section we present a theorem which gives a sufficient condition for a system of inequalities over $G$ to have a solution in $G$.
In the formulation we use Definition \ref{trans} and the notation given after it.  
In particular recall that if $w_{\mathcal{X}}(\bar{y})$ is oscillating but is not explicitly oscillating then 
\[ 
\hat{O}_{w}:=\bigcup _{V\in\mathcal{P} ^{os}} \Big( (V\setminus \mathsf{Fix}(G))\cap\bigcap _{i=0} ^{n_{V}-1} v_{V,0}^{-1}v
_{V,1}^{-1}\ldots v_{V,i}^{-1}\Big( \mathsf{supp}(v_{V,i+1})\Big)\Big) ,  
\] 
where for each $V\in\mathcal{P}^{os}$  
the word $w_{V}(\bar{y})$ is viewed in the form $\mathsf{u}_{V,n_{V}}v_{V,n_{V}}\ldots \mathsf{u}_{V,1}v_{V,1}$ with $v_{V,i}\in G\setminus\{ 1\}$
and  $\mathsf{u}_{V,i}$ depending only on variables, $i\leq n_{V}$ ( $v_{V,0}=1$). 
Also recall that in the case of explicitly oscillating $w_{\mathcal{X}} (\bar{y})$ we put $\hat{O}_w = O_w$.  

\begin{theorem} \label{uab}
Let $G$ act on a perfect metric space $(\mathcal{X},\rho )$ by
homeomorphisms. 
Let 
$\{ w_{1}(\bar{y}), w_{2}(\bar{y}), \ldots , w_{m}(\bar{y})\}$ 
be a set of reduced and non-constant words from  
$\mathbb{F} _{t}\ast G$ on $t$ variables  
$y_{1},\ldots ,	y_{t}$. 

If $G$ hereditarily separates $\mathcal{X}$ and every $w_{j}(\bar{y})$, $j\leq m$, is oscillating, then the set of inequalities 
$w_{1}(\bar{y})\neq 1, w_{2}(\bar{y})\neq 1,\ldots , w_{m}(\bar{y})\neq 1$ has a solution in $G$. 

Moreover, for any collection $\{ O_{j}\}$ such that $O_{j}$ is an open subset of the set $\hat{O}_{w_{j}}$, $j\leq m$, there is a solution $(g_{1},\ldots , g_{t})$ of this set of inequalities such that 
$\mathsf{supp}(g_{i})\subseteq\bigcup _{j=1}^{m} (\bigcup \mathcal{V}_{w_{j}}(O_{j}))$ for $1\leq i\leq t$.
\end{theorem}

\emph{Proof.} 
Wlog assume that when 
$w_{j}(\bar{y})\notin\mathbb{F} _{t}$, 
$w_{j}(\bar{y})=\mathsf{u}_{j,n_{j}}v_{j,n_{j}}\ldots \mathsf{u}_{j,1}v_{j,1}$ for $1\leq j\leq m$ with $v_{i}\in G\setminus\{ 1\}$. 
Fix the collection $\{ O_{j}\}$ from the statement of the theorem. 
For every $j \le m$ choose a point $o_j \in O_j$ so that 
$\mathcal{V} _{w_{j}}(\{ o_{j}\})\cap \mathcal{V} _{w_{j'}}(\{ o_{j'}\} )=\emptyset$ when $j,j' \in\{ 1,\ldots , m\}$ and $j\neq j'$. 
Using continuity of all $v_{j,s}$, for every $j\leq m$ we choose some open ball $B_{j}\subseteq O_{j}$ such that the
following conditions are satisfied: 
\begin{itemize} 
\item $o_j \in B_j$ for every $j\le m$ and  
\item  $(\bigcup \mathcal{V} _{w_{j}}(B_{j}))\cap (\bigcup \mathcal{V} _{w_{j'}}(B_{j'}))=\emptyset$ 
when $j,j' \in\{ 1,\ldots , m\}$ and $j\neq j'$. 
\end{itemize} 
Now we construct a sequence 
$(\bar{g} _{0}, \bar{g} _{1},\ldots ,\bar{g} _{m})$, where
$\bar{g} _{i}=(g_{i,1},\ldots , g_{i,t})$, such that for every $i$ with $1\leq i\leq m$ the following conditions are satisfied: 
\begin{itemize} 
\item $\mathsf{supp}(\bar{g _{i}})\subseteq\bigcup _{j=1}^{i}\mathcal{V}_{w_{j}}(B_{j})$, 
\item for all $j,\ell$ with $1\leq j\leq t$,  
$1\le \ell\le m$ and $\ell\neq i$ 
the restriction $g_{i,j}\upharpoonright_{B_{\ell}}$ coincides with 
$g_{i-1,j}\upharpoonright_{B_{\ell}}$, 
\item $\bar{g}_{i}$ is a solution of the set of inequalities $w_{1}(\bar{y})\neq 1,\ldots , w_{i}(\bar{y})\neq 1$.
\end{itemize}
Fix $\bar{g}_{0}:=(1,\ldots , 1)\in G^{t}$. 
Suppose that after $k-1$ steps the tuple 
$\bar{g}_{k-1}=(g_{k-1,1},\ldots , g_{k-1,t})$ is defined, $k\leq m$. 
At the $k$-th step we will modify the action of elements of this tuple on the ball $B_{k}$ so that $\bar{g}_{k}$ satisfies 
$w_{k}(\bar{g}_{k})\upharpoonright_{B_{k}}\neq \mathsf{id}\upharpoonright_{B_{k}}$ 
(i.e. $\bar{g}_{k}$ is a solution of the inequality 
$w_{k}(\bar{y})\neq 1$). 
For the $k$-th word $w_{k}(\bar{y})$ we consider two cases. 

\noindent
\textbf{Case 1.} 
$w_{k}(\bar{y})$ is explicitly oscillating. \\ 
Let us apply Theorem \ref{ab} to the word $w_{k}(\bar{y})$ and the set $B_{k}\subseteq O_{w_{k}}$. 
We obtain some solution 
$\bar{f}=(f_{1},\ldots , f_{t})\in G$ of the inequality 
$w_{k}(\bar{y}) \neq 1$, such that 
$\mathsf{supp}(f_{i})\subseteq \mathcal{V}^{> 0}_{w_{k}}(B_{k})$, $1\leq i\leq t$. 
Now for every $i\leq t$ we define $g_{k,i}:=f_{i}g_{k-1,i}$. 
Thus $\bar{g}_{k}$ is also a solution of 
$w_{k}(\bar{y})\neq 1$.
Since
\[ 
\bigcup _{i=1} ^{t} \mathsf{supp}(f_{i})\subseteq\mathcal{V}^{> 0}_{w_{k}}(B_{k})\mbox{ and } \bigcup_{i=1}^{t} \mathsf{supp}(g_{k-1,i})\subseteq\Big(\mathcal{X}\setminus\mathcal{V}_{w_{k}}(B_{k})\Big) , 
\]
the tuple $(g_{k,1},\ldots , g_{k,t})$ still is a solution of $w_{j}(\bar{y})\neq 1$ for $1\leq j\leq k-1$.  

\noindent 
\textbf{Case 2.} $w_{k}(\bar{y})$ is oscillating but is not explicitly oscillating. \\ 
Then applying $\mathsf{Transition}$ we obtain non-empty $\mathcal{P}^{os}$ corresponding to $w_{k}(\bar{y})$ and $B_k$. 
Thus there is some $U\in \mathcal{P}^{os}$, $U\subseteq B_k$,  and a word $w_{U}(\bar{y})$, which is derived from $w_{k}(\bar{y})$ by cancellations of constants and reductions, which is explicitly oscillating on $U$. 
By Theorem \ref{ab} there is some $\bar{f}$ such that for every $i\leq t$, $\mathsf{supp}(f_{i})\subseteq \mathcal{V}_{w
_{k}}(B_{k})$, for every $j\leq n_{k}$, $f_{i}$ stabilizes $v_{j}\ldots v_{1}(U)$ setwise and 
$w_{U}(\bar{f})(p)\neq p$ for some point $p\in U$. 
Thus by Lemma \ref{ii} we have $w_{k}(\bar{f})\neq 1$. 

 Now similarly as above for every $i\leq t$ we define $g_{k,i}:=f_{i}g_{k-1,i}$. 
This gives a solution of all inequalities $w_{j}(\bar{y})\neq 1$ for $j\leq k$. 

 Thus, after $m$ steps of the algorithm we obtain a tuple $\bar{g} _{m}$, which is the solution of
the system $w_{1}(\bar{y})\neq 1,\ldots , w_{m}(\bar{y})\neq 1$. Moreover, for any $1\leq i\leq t$, 
$\mathsf{supp}(g_{i})\subseteq\bigcup _{j=1} ^{m}\mathcal{V}_{w_{j}}(O_{j})$. $\square$ 

\bigskip 

\begin{remark} \label{diameter2} 
{\em The conclusion of Theorem \ref{uab} can be extended by the following additional statement:  
\begin{itemize} 
\item for any $\varepsilon$ the solution $\bar{g}$ can be chosen so that it additionally satisfies the inequality 
$d(x,g_i(x))\le \varepsilon$ for all $i\le t$ and $x\in \mathcal{X}$.  
\end{itemize}  
To see this it is enough to add the following argument at each step of the inductive procedure of the proof of Theorem \ref{uab}. 
Using Remark \ref{diameter1} at each step of the induction we choose the corresponding solution $\bar{f}$ of  
$w_k (\bar{y}) \not=1$ with the additional property that 
$d(f_i (x) ,x)\le \varepsilon$ for all $i\le t$ and 
$x\in \mathcal{X}$. 
Then one easily sees by inspection of the proof that the choice of the family $\{ B_j \, | \, 1\le j\le m \}$ guarantees that $d(g_i (x) , x)\le \varepsilon$ for all $j\le t$ and $x\in \mathcal{X}$. } 
\end{remark} 

\begin{remark} 
{\em One can wonder if in order to solve the system of inequations given in Theorem \ref{uab} it is possible to apply the trick of commutators used in the proof of Proposition \ref{e_c}. 
We should mention here that the situation of Proposition \ref{e_c} is very special, where the commutator of oscillating words is again an oscillating word. 
This is not true in general and cannot be applied in Theorem \ref{uab}. } 
\end{remark}

\section{Actions of topological groups} 

\subsection{Polish $G$-spaces}

A {\em Polish space (group)} is a separable, completely
metrizable topological space (group).
The corresponding metric extends to tuples by 
\[ 
d((x_1 ,...,x_m ), (y_1 ,...,y_m ))= \mathsf{max} (d(x_1 ,y_1 ),...,d(x_m ,y_m )). 
\]  
Let $(\mathcal{X},d)$ be a Polish space and $\mathsf{Iso}(\mathcal{X})$ 
be the corresponding isometry group 
endowed with the pointwise convergence topology. 
Then $\mathsf{Iso} (\mathcal{X})$ is a Polish group. 
A compatible left-invariant metric can be obtained as follows.  
Fix a countable dense set $S=\{ s_i : i\in \{ 1,2,...\} \} \subseteq \mathcal{X}$.  
Define for two isometries $\alpha$ 
and $\beta$ of $\mathcal{X}$ 
\[  
\rho_{S} (\alpha ,\beta )= \sum_{i=1}^{\infty} 2^{-i} \mathsf{min}(1, d(\alpha (s_i ),\beta (s_i ))) .
\]  
Let $G$ be a closed subgroups of $\mathsf{Iso}(\mathcal{X})$.  
We fix a dense countable set $\Upsilon \subset G$ and  
distinguish a base of $G^t$ consisting of all sets of the following form. 
Let $\bar{s}_1 ,\ldots ,\bar{s}_t$ be a sequence of tuples from $S$, $q\in \mathbb{Q}$ and $h_1 ,\ldots , h_t$ be a sequence from $\Upsilon$. 
Define 
\[ 
N^q (\bar{s}_1 ,h_1 ,\ldots ,\bar{s}_t, h_t ) = \{ (g_1 , \ldots , g_t ) \, | \, d(g_i (\bar{s}_i ),h_i (\bar{s}_i ) <q \, , 1 \le i \le t  \}.   
\] 
The family of all 
$N^q (\bar{s}_1 ,h_1 ,\ldots ,\bar{s}_t, h_t )$ 
forms a base of the topology of $G^t$. 
This material can be found in any textbook from  descriptive set theory, see for example \cite{kechris}. 

The following theorem is related to Theorem \ref{uab} and Remark \ref{diameter2}. 

\begin{theorem} \label{density} 
Let $G$ be a closed subgroup of the isometry group 
$\mathsf{Iso} (\mathcal{X})$ of a perfect Polish space 
$\mathcal{X}$. 
Assume that the action of $G$ is hereditarily separating on $\mathcal{X}$.  
Then for any oscillating word $w(\bar{y})$ from  
$\mathbb{F} _{t}\ast G$ on $t$ variables, 
$y_{1},\ldots ,	y_{t}$, the set  
$\{ \bar{g} \, | \, w(\bar{g})\neq 1 \}$ is dense in $G^t$. 
\end{theorem} 

{\em Proof.} 
Let us fix an inequality $w(\bar{y}) \not= 1$ over $G$ 
and assume that $w(\bar{y})$ is explicitly oscillating. 
If $w(\bar{y}) \not\in \mathbb{F}_t$ we may assume 
that it is in the form (1.2). 
Let $N^q (\bar{s}_1 ,h_1 ,\ldots ,\bar{s}_t, h_t )$ be 
a basic open set defined before the theorem. 
Let $p$ be an element of $\mathcal{X}$ 
which does not occur in $\bar{s}_1 ,\ldots ,\bar{s}_t$. 
Define  
\[ 
U^q (\bar{s}_1 ,h_1 ,\ldots ,\bar{s}_t, h_t ,p, w(\bar{y})) 
\] 
to be the set of tuples from 
$N^q (\bar{s}_1 ,h_1 ,\ldots ,\bar{s}_t, h_t )$ 
which are distinctive for $w(\bar{y})$ and $p$. 
It is easy to see that $U^q$ is open in $N^q$. 
For example one can repeat the argument given in p.530 of \cite{A} (extended by application of homeomorphisms $v_j$). 

Assume that 
$N^q (\bar{s}_1 ,h_1 ,\ldots ,\bar{s}_t, h_t )\not=\emptyset$. 
We will prove that if  
$p\in O_w \setminus \bigcup \mathcal{V}^{-1}_{w}(\bigcup \mathcal{V}_w (\overline{O}_w \setminus O_w ))$  
then the set 
$U^q (\bar{s}_1 ,h_1 ,\ldots ,\bar{s}_t, h_t ,p, w(\bar{y}))$   
is not empty. 
Since this can be applied to every $q$ and every tuple 
$\bar{s}'_1 ,h'_1 ,\ldots ,\bar{s}'_t, h'_t$ such that $\bar{s}_i\subseteq \bar{s}'_i$ and $h'_i$ agrees with $h_i$ on $\bar{s}_i$, $1\le i \le t$,   
we would obtain that 
$U^q (\bar{s}_1 ,h_1 ,\ldots ,\bar{s}_t, h_t ,p, w(\bar{y}))$
is dense in $N^q (\bar{s}_1 ,h_1 ,\ldots ,\bar{s}_t, h_t )$. 

The case $|w(\bar{y})|=0$ is degenerate; we have $N^q = U^q$.  
We now apply some ideas from the proof of Theorem \ref{ab}.  
At $k$-th step of induction we will show that: 
\begin{itemize} 
\item There is a tuple 
\[ 
(g_{1},\ldots , g_{t})
\in U^q (\bar{s}_1 ,h_1 ,\ldots ,\bar{s}_t, h_t ,p, (w)_k (\bar{y})),   
\]   
\item In the condition above we can choose $\bar{g}$ so that for all $i$ with $1\leq i\leq t$, 
$\mathsf{supp}(g_{i})\subseteq (\bigcup \mathcal{V}^{>0}_w (O_w ))\setminus (\bigcup \mathcal{V}^{-1}_{[w]_{k}} (\bigcup \mathcal{V}_w (\overline{O_w }\setminus O_w )))$ 
and each member of 
$\mathcal{V}_{w}(O_w )$ is $g_i$-invariant. 
\end{itemize}  
We start by fixing some  
$( g^{\circ}_1, g^{\circ}_2, \ldots , g^{\circ}_t )\in N^q (\bar{s}_1 ,h_1 ,\ldots ,\bar{s}_t, h_t )$ and $q'<q$ such that 
$d(g^{\circ}_i (\bar{s}_i ),h_i (\bar{s}_i )) <q'$ for all 
$i$ with $1 \le i \le t$.  
We will assume that $w(\bar{y}) \notin \mathbb{F}_t$. 
When $w(\bar{y}) \in \mathbb{F}_t$ the argument below works.   In fact after removing the corresponding $v_i$ it becomes easier. 
In particular, we assume that $(w)_{1}(\bar{y})$ is of the form 
$y_{j}^{\pm 1}v_{1}$ for some $1\leq j\leq t$.  
When $w(\bar{y})\not\in\mathbb{F} _{t}$ then 
$p'\neq v_{1}(p')$ for all  $p'\in O_w$. 
Thus according to the assumptions on $p$, for a non-trivial $v_1$ the inequality $p\neq v_{1}(p)$ is satisfied. 

Wlog suppose $(w)_1 (\bar{y}) =y_{1}v_{1}$. 
If $g^{\circ}_{1}v_{1}(p) \notin \{ p, v_1 (p)\}$ then 
\[ 
( g^{\circ}_1, g^{\circ}_2, \ldots , g^{\circ}_t )\in U^q (\bar{s}_1 ,h_1 ,\ldots ,\bar{s}_t, h_t ,p, (w)_1 (\bar{y})).  
\] 
In the contrary case using the argument of step $k=1$ of the proof of Theorem \ref{ab} and Remark \ref{diameter1} find 
\[ 
( g'_1, g'_2, \ldots , g'_t ) \in 
U^{q-q'} (\bar{s}_1 ,\mathsf{id} ,\ldots ,\bar{s}_t, \mathsf{id} ,p, (w)_1 (\bar{y}))
\] 
such that 
\[ 
(w)_1 (\bar{g}')(p) \notin \Big\{ p, v_1 (p),(g^{\circ}_1)^{-1} (p), (g^{\circ}_1)^{-1} (v_1 (p))\Big\}. 
\] 
Then 
\[ 
( g^{\circ}_1g'_1 , g^{\circ}_2 g'_2 , \ldots , g^{\circ}_t g'_t )\in U^q (\bar{s}_1 ,h_1 ,\ldots ,\bar{s}_t, h_t ,p, (w)_1 (\bar{y})),   
\] 
i.e. we may take $g_i = g^{\circ}_i g'_i$ , $1\le i\le n$.  
It is worth noting that instead of $q - q'$ one could take arbitrary small $q''$. 
In particular $( g^{\circ}_1g'_1 , g^{\circ}_2 g'_2 , \ldots , g^{\circ}_t g'_t )$ can be chosen arbitrary close to $( g^{\circ}_1, g^{\circ}_2 , \ldots , g^{\circ}_t )$. 

If $|w(\bar{y})|>2$ then for the second step we additionally need that \\ 
$g_1 ((v_{1}(p))\notin \bar{s}_1 \cup \ldots \cup \bar{s}_t$. 
This can be arranged at the first step by the stronger demand that 
\[ 
g'_1 (v_1 (p)) \notin \Big\{ p, v_1 (p),(g^{\circ}_1)^{-1} (p), (g^{\circ}_1)^{-1} (v_1 (p))\Big\}\cup (g^{\circ}_1)^{-1} (\bar{s}_1) \cup \ldots \cup (g^{\circ}_1)^{-1} (\bar{s}_t). 
\] 
In fact our argument shows that the set 

$\Big\{ (g_1, g_2, \ldots , g_t )\in U^q (\bar{s}_1 ,h_1 ,\ldots ,\bar{s}_t, h_t ,p, (w)_1 (\bar{y})) \, | \, $ 
 
\hspace{5cm} $(w)_1 (\bar{g})(p))\notin\{ p, v_{1}(p)\} \cup \bar{s}_1 \cup \ldots \cup \bar{s}_t \} \Big\}$ \\    
is dense in $N^q (\bar{s}_1 ,h_1 ,\ldots ,\bar{s}_t, h_t)$. 
  
For a fixed $\bar{g}$ from this part of 
$U^q( \bar{s}_1 ,h_1 ,\ldots ,\bar{s}_t, h_t ,p, (w)_1 (\bar{y}))$  take $q_1 < q$ with 
$\bar{g} \in U^{q_1}( \bar{s}_1 ,h_1 ,\ldots ,\bar{s}_t, h_t ,p, (w)_1 (\bar{y}))$.  
At  Step 2 we look for an element of    
\[ 
U^{q-q_1}(v_1 (p)\bar{s}_1 ,g_1 ,\ldots ,\bar{s}_t, g_t ,p, (w)_2 (\bar{y})) \subset N^{q} (\bar{s}_1 ,h_1 ,\ldots ,\bar{s}_t, h_t ).   
\]  
Skipping further considerations of this step we jump to Step $k$. 

After Step $k-1$ we have some $q_{k-1} < q$ and some tuple $\bar{g}$ in 
\[ 
U^{q_{k-1}}(\bar{s}'_1 ,h'_1 ,\ldots ,\bar{s}'_t, h'_t ,p, (w)_{k-1} (\bar{y})) \subseteq N^{q} (\bar{s}_1 ,h_1 ,\ldots ,\bar{s}_t, h_t ).    
\] 
where $\bar{s}_i \subseteq \bar{s}'_i$, $i\le t$ and each $h'_i$ agrees with $h_i$ on $\bar{s}_i$. 
We look for a tuple from   
\[ 
U^{q-q_{k-1}}(\bar{s}'_1 ,h''_1 ,\ldots ,\bar{s}'_t, h''_t ,p, (w)_{k} (\bar{y}))   
\]   
where each $h''_i$ agrees with $g_i$ on $\bar{s}'_i$. 
Here we again have two cases. 

\noindent
{\em Case 1.} $(w)_{k}(\bar{y}) =u_{d,s+1}u_{d,s}\ldots u_{d,1}v_{d}\ldots u_{2,1}v_{2}u_{1,\ell_{1}}\ldots u_{1,1}v_{1}$, 

\hspace{7cm} where $k-1=L_{d-1}+s$, $s\geq 1$.\\ 
Our argument is a slight modification of the corresponding place in the proof of Theorem \ref{ab}. 
We preserve the notation of this proof.  
In particular for a fixed $\bar{g}$ as above we take the corresponding $p_{i,\bar{g}}$, $i\le k$. \\ 
If
\[ 
p_{k,\bar{g}}\notin\Big\{ p_{i,\bar{g}}\ \Big|\ 0\leq i\leq k-1\Big\}\cup\Big\{ v_{1}(p_{0,\bar{g}}),\ldots
	, v_{d}(p_{L_{d-1},\bar{g}})\Big\} , 
\] 
then we have found an acceptable tuple $\bar{g}$. 
Let us assume that $p_{k,\bar{g}}=p_{m,\bar{g}}$ for some
$0\leq m<k$ or $p_{k,\bar{g}}=v_{m+1}(p_{L_{m},\bar{g}})$ for some $0\leq m<d-1$. 

As in the proof of Theorem \ref{ab} we may assume that 
$u_{d,s+1}=y_{j}$ and put again 
\[ 
Y:=\Big\{ p_{i,\bar{g}}\ \Big|\ 0\leq i\leq k-2\Big\}\cup\Big\{ v_{1}(p_{0,\bar{g}}),\ldots , v_{d}(p_{L_{d-1},
	\bar{g}})\Big\} . 
\] 
The neighborhood $O\subseteq v_{d}\ldots v_{1}(O_w )$ of the point 
$p_{k-1,\bar{g}}\in v_{d}\ldots v_{1}(O_w )$ is chosen as before. 
We can take it sufficiently small so that when we define  
$f\in \mathsf{stab}_{G}((X\setminus O)\cup \bigcup \mathcal{V}^{-1}_{[w]_{k-1}} (\bigcup \mathcal{V}_w (\overline{O_w}\setminus O_w))\cup Y)$, 
we have that 
$\mathsf{max}_{1\le i \le t}\{ d(\bar{s}'_i , f(\bar{s}'_i)) \} \le q -  q_{k-1}$. 
Since we want 
$p_{k, \bar{g}} \notin Y\cup \bar{s}_1 \cup \ldots \cup \bar{s}_t$, 
we arrange $f$ taking $p_{k-1,\bar{g}}$ outside $Z$, where 
\[ 
Z:=\Big\{ g_{j}^{-1}(p_{i,\bar{g}})\ \Big| \ 0\leq i\leq k-1\Big\}\cup\Big\{ g_{j}^{-1}(v_{i+1}(p_{L_{i},\bar{g}}))\
\Big|\ 0\leq i\leq d-1\Big\} \cup 
\] 
\[
\hspace{8cm} 	(g_j)^{-1} (\bar{s}'_1) \cup \ldots \cup (g_j)^{-1} (\bar{s}'_t).  
\]
We apply hereditary separation at this point. 
Replacing $g_{j}$ by $g_{j}f$ we obtain a corrected tuple $\bar{g}$.  
By the choice of $q_{k-1}$ and $O$ this tuple represents  
\[ 
U^{q-q_{k-1}}(\bar{s}'_1 ,h''_1 ,\ldots ,\bar{s}'_t, h''_t ,p, (w)_{k} (\bar{y}))   
\subset U^{q}(\bar{s}'_1 ,h'_1 ,\ldots ,\bar{s}'_t, h'_t ,p, (w)_k (\bar{y})). 
\] 
We finish the proof of Case 1 as in Theorem \ref{ab}. 

\noindent 
{\em Case 2.} 
$(w)_{k}=u_{d+1,1}v_{d+1}u_{d,s}\ldots u_{d,1}v_{d}\ldots u_{2,1}v_{2}u_{1, \ell_{1}}\ldots u_{1,1}v_{1}$, 
where $k=L_{d}+1$.
This corresponds to Case 2 of the proof of Theorem \ref{ab} and it should be treated in the same fashion. 
The rest of the argument for explicitly oscillating $w(\bar{y})$ is clear. 

Let us consider the case when $w(\bar{y})$ is oscillating 
but is not explicitly oscillating.  
Then applying $\mathsf{Transition}$ we obtain non-empty $\mathcal{P}^{os}$ corresponding to $w(\bar{y})$. 
Thus there is some word $w_{V}(\bar{y})$ derived
from $w(\bar{y})$ by cancellations of constants and reductions, which is explicitly oscillating. 
Then fixing $N^q (\bar{x}_1 ,h_1 ,\ldots ,\bar{x}_t, h_t )$ we repeat the proof above replacing $w(\bar{y})$ by $w_V (\bar{y})$. 
Take a point $p\in V$ according the requirements of that proof. 
Then we see that  
$U^q (\bar{x}_1 ,h_1 ,\ldots ,\bar{x}_t, h_t ,p, w_V (\bar{y}))$ 
is dense in  
$N^q (\bar{x}_1 ,h_1 ,\ldots ,\bar{x}_t, h_t )$. 
In particular, so is the set  
$\{ (g_1, g_2, \ldots , g_t )\in N^q \, | \, w_V (\bar{g}) \not= 1 \}$.    
Thus by Lemma \ref{ii} we have that $\{ (g_1, g_2, \ldots , g_t )\in N^q \, | \, w(\bar{g}) \not= 1 \}$ is dense in  
$N^q (\bar{x}_1 ,h_1 ,\ldots ,\bar{x}_t, h_t )$. 
$\Box$ 

\subsection{Topological actions of locally compact groups} 

A locally compact group $G$ carries a left $G$-invariant Haar measure $\mu$. 
An action of $G$ on a set $X$ is called {\em topological} if each point stabilizer $\mathsf{stab}_G(x)$, $x\in X$, is closed and of measure $0$.  
An action is called {\em strongly topological} if for any finite $Y\subseteq X$ each point stabilizer $\mathsf{stab}_G(Y \cup \{ x\})$,  $x\in X\setminus Y$, is closed in $\mathsf{stab}_G(Y )$ and has measure $0$ under the Haar measure for $\mathsf{stab}_G(Y )$.  
Under the assumption that all point-stabilizers are non-trivial this condition strengthens separating actions. 

In this section we concentrate on strongly topological actions of locally compact groups. 
We will assume that $X$ is a Polish space, say $(\mathcal{X},d)$,  and $G$ acts on $\mathcal{X}$ isometrically and strongly topologically. 
We preserve notation of Section 4.1. 

If $H$ is a closed subgroup of $G$, it has a left Haar measure $\lambda$. 
If $K$ is a closed subgroup of $H$ of zero measure (and of infinite index), then the space $H/K$ inherits a natural measure $\bar{\lambda}$ which is equivariant, 
i.e. $H$-translates of $\bar{\lambda}$-negligible sets are $\bar{\lambda}$-negligible. 
As $K$ is a closed subgroup of $H$ it also has a left Haar measure, say $\mu_K$. 

The following lemma is Lemma 3.1 from \cite{A}. 

\begin{lemma} \label{A3_1}
Let $H_0$ be a $\lambda$-measurable subset of $H$ such that for almost all $K$-cosets $D\subseteq H$ we have 
$\mu_K (D\cap H_0 )=0$. 
Then $\lambda (H_0 )=0$.  
\end{lemma} 

\noindent
We will apply this lemma to left cosets of $H$ instead of $H$. 
These cosets are considered under an obvious extension of the  measure of $H$.  
Then $H_0$ from the formulation will be a subset of some $N^q (s, \mathsf{id})$, $s\in S$. 
Typically $K$ from the lemma arises as 
the stabilizer of some $p\in \mathcal{X}$. 
Under our assumptions $\mu (\mathsf{stab}_G(p))=0$. 

The following theorem is a version of Theorem 1.5 of 
\cite{A}. 
It looks slightly technical, but in Corollary \ref{clc} 
we give a short and natural formulation. 

\begin{theorem} \label{lc} 
Let $G$ be a locally compact topological group acting strongly topologically on a Polish space $(\mathcal{X},d)$ by isometries. 
Assume that for every finite $Y\subset \mathcal{X}$ the stabilizer $\mathsf{stab}_G(Y)$ is not trivial. 
Let $s\in S$, $q\in \mathbb{Q}$. 
\begin{itemize}  
\item Let $w(\bar{y})$ be an explicitly oscillating word such that $O_w$ contains \\ 
$B_{(|w(\bar{y})|+1)q} (s)$ and each $B_{\ell q} (s)$ with $1\le \ell \le |w(\bar{y})|+1$, is infinite and invariant with respect to each constant $v_i$ of the word $w(\bar{y})$.      
\item Let $\bar{\gamma}$ be the random $t$-tuple in $N^q (s,\mathsf{id})$. 
\end{itemize} 
Then almost surely $\bar{\gamma}$ satisfies the mixed inequality $w(\bar{y}) \not= 1$.  
\end{theorem} 

{\em Proof.} 
We adapt the proof of Theorem 1.5 from \cite{A}. 
Let $\bar{x}_1 ,\ldots ,\bar{x}_t$ be a sequence of tuples from  
$\mathcal{X}$ and $h_1 ,\ldots , h_t$ be a sequence from 
$N^q (s,\mathsf{id})$. 
Define 
\[ 
A(\bar{x}_1 ,h_1 ,\ldots ,\bar{x}_t, h_t ) = \{ (g_1 , \ldots , g_t ) \, | \, 
g_i (\bar{x}_i )= h_i (\bar{x}_i ) \, , 1 \le i \le t  \}.   
\]
We view this set as follows. 
Let $G_0$ be the direct product of stabilizers 
$\mathsf{stab}_G(\bar{x}_1) \times \ldots \times \mathsf{stab}_G(\bar{x}_t)$.  
Then $A(\bar{x}_1 ,h_1 ,\ldots ,\bar{x}_t, h_t )$ is a right coset of $G_0$. 
It inherits the natural topology and measure from $G_0$.  

In order to introduce the next notation used in the proof let $w(\bar{y})$ be an explicitly oscillating word and $\ell$ be a natural number. 
If $w(\bar{y}) \not\in \mathbb{F}_t$ we may assume that it is in the form (1.2). 
Now let $p$ be an element of $B_{\ell q} (s)$ such that neither $p$ nor $v_1 (p)$ occurs in $\bar{x}_1 ,\ldots ,\bar{x}_t$. 
Put 
\[ 
U(\bar{x}_1 ,h_1 ,\ldots ,\bar{x}_t, h_t ,p, w(\bar{y})) 
\] 
to be the set of tuples from 
$A(\bar{x}_1 ,h_1 ,\ldots ,\bar{x}_t, h_t )$ which are distinctive for $w(\bar{y})$ and $p$ and the elements 
\[ 
p=p_{0,\bar{g}}, v_{1}(p_{0,\bar{g}}),\ldots , p_{l_{1},\bar{g}}, v_{2}(p_{l_{1},\bar{g}}),\ldots , p_{L_n , \bar{g}}
\] 
do not occur in $\bar{x}_1 , \ldots ,\bar{x}_t$. 
Then $U$ is open in $A$. 
This situation is similar to the one appeared in the proof of Theorem \ref{density}. 
Again in order to show openess  we can repeat the argument given on p. 530 of \cite{A} (extended by application of homeomorphisms $v_j$). 

Applying induction on $k=|w(\bar{y})|$ we will prove that 
\begin{itemize} 
\item 
if $B_{(k+\ell)q}(s) \subset O_{w}$ , 
every $B_{mq} (s)$ with $1 \le m \le  k+ \ell +1$, is infinite and invariant with respect to each constant $v_i$ of the word $w(\bar{y})$, 
\[ 
\bar{x}_1 \cup \ldots \cup \bar{x}_t \subset B_{\ell q}(s) \mbox{ and }  
\{ p, v_1 (p)\} \subset  
B_{\ell q}(s) \setminus ( \bar{x}_1 \cup \ldots \cup \bar{x}_t )  
\]   
then the set 
\[ 
 U(\bar{x}_1 ,h_1 ,\ldots ,\bar{x}_t, h_t ,p, w(\bar{y})) 
\]    
is almost surely in 
$(N^q (s, \mathsf{id}))^t \cap A(\bar{x}_1 ,h_1 ,\ldots ,\bar{x}_t, h_t )$, i.e. $\mu ((N^q (s, \mathsf{id})^t \cap A)\setminus U)=0$. 
\end{itemize} 
The case $|w(\bar{y})|=0$ is degenerate and obvious: we have $A = U$.  
When $w(\bar{y}) \in \mathbb{F}_t$ the argument coincides with the corresponding one from the proof of Theorem 1.5 in \cite{A}.  
Basically this is a simplified version the argument below for  
$w(\bar{y}) \not\in \mathbb{F}_t$. 
We will also apply some ideas from the proof of Theorem \ref{ab}. 

Consider Step $k>0$.
We assume that the statement 
$\mu (((N^q (s,\mathsf{id}))^t\cap A)\setminus U)=0$  holds for all words of length $k' < k$,  all tuples of parameters 
$\bar{x}'_1 ,h'_1 ,\ldots ,\bar{x}'_t, h'_t$ with 
$\bar{x}'_1 ,\ldots ,\bar{x}'_t \in B_{\ell q} (s)$ and all $p$ taken as it was described above.  
We will assume that $(w)_{1}$ is of the form $y_{i}^{\pm 1}v_{1}$ for some $1\leq i\leq t$.  
Our argument also works in the case $v_1  =1$ (which is now allowed). 
Note that in the case $v_1 \not=1$ according to the assumptions on $p$ we have $p\neq v_1 (p)$ and $v_1 (p)\in B_{\ell q}(s)$. 

Wlog suppose $(w)_1 =y_{i}v_{1}$. 
The set 
\[ 
\Big\{ (g_1, g_2, \ldots , g_t )\in A \, | \, g_i (v_1 (p)) \in 
\bar{x}_1 \cup h_1 (\bar{x}_1 )\ldots \cup \bar{x}_t \cup h_t(\bar{x}_t ) \cup \{ p, v_1 (p)\}\Big\} 
\]  
is a finite union of cosets of the direct product  
\[ 
S_{v_1 (p)} = \mathsf{stab}_G(\bar{x}_1) \times \ldots \times \mathsf{stab}_G(\{ v_1 (p) \}\cup \bar{x}_i) \times \ldots \times \mathsf{stab}_G(\bar{x}_t).
\] 
Since $G$ acts strongly topologically on $\mathcal{X}$, we deduce that for almost all \\ 
$(g_1, g_2, \ldots , g_t )\in (N^q (s,\mathsf{id}))^t\cap A$ we have 
\[ 
g_i (v_1(p)) \not\in 
\bar{x}_1 \cup h_1(\bar{x}_1) \cup  \ldots \cup \bar{x}_t \cup h_t(\bar{x}_t )\cup \{ p,v_1 (p)\}. 
\] 
For such a tuple and for a fixed $p'$ of the form $g_i ((v_{1}(p))$ let $h'_i$ be any element of $N^q (s, \mathsf{id})$ mapping $\{ v_1 (p) \}\cup \bar{x}_i$ to $\{ p'\} \cup h(\bar{x}_i )$. 
Note that $p'\in B_{(\ell +1)q}(s)$. 

Consider  
\[ 
U(\bar{x}_1,h_1 ,\ldots ,\{ v_1 (p)\} \cup\bar{x}_i ,h'_i ,\ldots ,\bar{x}_t, h_t ,p', [w]_1 (\bar{y}))
\] 
in $(N^q (s, \mathsf{id}))^t \cap A(\bar{x}_1, h_1 ,\ldots ,\{ v_1 (p)\} \cup \bar{x}_i ,h'_i ,\ldots ,\bar{x}_t, h_t )$.    
Since $B_{(\ell + 1 + |[w]_1|)q}(s) \subseteq O_{[w]_1}$ (by the inductive assumptions and assumptions on constants in $w(\bar{y})$) and 
\[ 
p' \in B_{(\ell +1)q}(s) \setminus 
(\{ p, v_1 (p) \} \cup \bar{x}_1 \cup \ldots \cup \bar{x}_t ) 
\] 
applying induction we see that  
$U(\bar{x}_1 ,h_1 ,\ldots , \{ v_1 (p)\} \cup \bar{x}_i ,h'_i ,\ldots ,\bar{x}_t, h_t ,p', [w]_1 (\bar{y}))$ 
is almost surely in 
$(N^q (s, \mathsf{id}))^t \cap A(\bar{x}_1 ,h_1 ,\ldots ,\{ v_1 (p)\} \cup \bar{x}_i ,h'_i ,\ldots ,\bar{x}_t, h_t )$. 
In particular $U(\bar{x}_1 ,h_1 ,\ldots ,\bar{x}_i ,h_i ,\ldots ,\bar{x}_t, h_t ,p, w (\bar{y}))$  is almost surely in 
$(N^q (s, \mathsf{id}))^t \cap A(\bar{x}_1 ,h_1 ,\ldots ,\{ v_1 (p)\} \cup \bar{x}_i ,h'_i ,\ldots ,\bar{x}_t, h_t )$. 
Using Lemma \ref{A3_1} in the situation when $K= \mathsf{stab}_G(\bar{x}_1 ,\ldots ,\{ v_1 (p)\} \cup \bar{x}_i, \ldots, \bar{x}_t)$ 
we conclude that \\ 
$U(\bar{x}_1 ,h_1 ,\ldots ,\bar{x}_t, h_t ,p, w(\bar{y}))$ is almost surely in  
$(N^q (s, \mathsf{id}))^t \cap A(\bar{x}_1 ,h_1 ,\ldots ,\bar{x}_t, h_t )$. 

Let us fix a mixed inequality $w(\bar{y}) \not= 1$ and assume that $w(\bar{y})$ is explicitly oscillating and satisfies the conditions of the formulation of the theorem. 
Take $p\in B_q(s)$. 
Apply the claim proved by induction to $A=G^t$ and 
$U(p, w(\bar{y}))$ where $\ell = 1$. 
We obtain that $\mu ((N^q (s, \mathsf{id}))^t \setminus U)=0$. 
$\Box$

\bigskip

The following statement is proved by the same proof as in Theorem \ref{lc}. 
One only has to replace balls $B_{\ell q}(s)$ by $\mathcal{X}$ and $N^q (s, \mathsf{id})$ by $G$.

\begin{corollary} \label{clc} 
Let $G$ be a locally compact topological group acting strongly topologically on a Polish space $(\mathcal{X},d)$ by isometries. 
Assume that for every finite $Y\subset \mathcal{X}$ the stabilizer $\mathsf{stab}_G(Y)$ is not trivial. 
Let $\bar{\gamma}$ be the random $t$-tuple in $G$ and $w(\bar{y})$ be an   explicitly oscillating word such that $O_w = \mathcal{X}$. 
Then almost surely $\bar{\gamma}$ satisfies the mixed inequality $w(\bar{y}) \not= 1$.  
\end{corollary}

\begin{corollary} \label{c}
Let $G$ be a compact topological group acting topologically on a Polish space $(\mathcal{X},d)$ by isometries. 
Assume that the action is separating. 
Let $\bar{\gamma}$ be the random $t$-tuple from $G$ and $w(\bar{y})$ be  an explicitly oscillating word such that $O_w = \mathcal{X}$. 
Then almost surely $\bar{\gamma}$ satisfies the mixed inequality $w(\bar{y}) \not= 1$.  
\end{corollary} 

{\em Proof.} 
To apply Theorem \ref{lc} we only have to verify that the action of $G$ on $\mathcal{X}$ is strongly topological. 
This verification coincides with the argument of Theorem 1.3 of \cite{A}. 
$\Box$ 

\bigskip 

\begin{corollary} 
Let $G$ be a profinite weakly branch group.  
Let $G$ act on a rooted tree $T$ spherically transitively such that the rigid vertex stabilizers are non-trivial. 
Let $\mathcal{X}$ be the boundary of $T$. 
Let $\bar{\gamma}$ be the random $t$-tuple from $G$ and $w(\bar{y})$ be an  explicitly oscillating word such that $O_w = \mathcal{X}$.   
Then almost surely $\bar{\gamma}$ satisfies the mixed inequality $w(\bar{y}) \not= 1$.  
\end{corollary}

{\em Proof.} 
Since $G$ is closed in the profinite topology of $\mathsf{Aut}(T)$, each stabilizer $\mathsf{stab}_G (x)$, $x\in \mathcal{X}$, is closed. 
By Example \ref{grig} $G$ hereditarily separates $\mathcal{X}$. 
By Corollary \ref{c} we have the conclusion of the corollary. 
$\Box$

\subsection{Topological actions of automorphism groups} 

Here we will consider the situation of Section 4.2 in the context of paper \cite{EGMM}. 
A word over a group $G$ with a single variable is of the following form: 

\bigskip 

$w(y) = h_k y^{m_k} \cdot \ldots \cdot h_1 y^{m_1}$ , where  
$h_1, \ldots , h_k \in G$  
\hspace{3cm} \, (4.1) 

\hspace{7cm} 
and $m_1, \ldots , m_k 
\in \mathbb{Z}\setminus \{ 0 \}.$   

\noindent  
We admit the possibility that $h_k =1$. 
According to Remark 5.1 of \cite{HO} in order to show that $G$ is MIF it suffices to prove that there are no laws of the form $w(y) = 1$ where $w(y)$ is as above. 
Section 6 of \cite{EGMM} gives a method of analysis of such words in the case of the automorphism group of some standard continuous structures. 
In fact the authors detect a certain (model-theoretic) property of these structures which guarantees that any $w(y)$ of the form $(4.1)$ is not a law. 
Adapting this property in the general situation of $G$-spaces we arrive at the following definition. 

\begin{definition} 
Assume that $G$ acts on an infinite set $X$ by permutations. 
We say that $(G,X)$ is  discerning if for every finite $A\subset X$, every non-trivial $\mathsf{stab}_G(A)$-orbit $O$ and every non-trivial $g\in G$ the intersection 
$\mathsf{supp}(g) \cap O$ is not empty.  
\end{definition}   
Lemma 6.13 of \cite{EGMM} states: 
\begin{quote} 
Assume that $(G,X)$ is discerning and $g_1, g_2 \in G$. 
Then for every finite $A \subset X$ and every $a\in X\setminus A$ there is $g\in g_2 \cdot \mathsf{stab}_G (A)$  such that 
$\{ g(a), g_1 g (a)\} \subset \mathsf{supp}(g_1 )\setminus (A \cup \{ a \})$.  
\end{quote} 
In \cite{EGMM} this property is applied to constants $h_i$ of the word $(4.1)$.  
These $h_i$ are viewed as $g_1$ in the formulation. 

The assumption that $G$ is a locally compact group with a strongly topological action on a space $\mathcal{X}$ already implies a statement of this kind: for every finite $A \subset \mathcal{X}$ and every $a\in \mathcal{X}\setminus A$ almost surely all elements $g$ of the coset $g_2 \cdot \mathsf{stab}_G (A)$ have the property that 
$\{ g(a), g_1 g (a)\} \subset \mathcal{X}\setminus (A \cup \{ a \})$.   
Indeed, since $\mathsf{stab}_G (A\cup \{ a \})$ is of measure 0 in $\mathsf{stab}_G (A)$ we see that so is the set 
\[ 
\{ g\, | \, g(a) \in (A \cup \{ a \}) \mbox{ or } g(a) \in 
 g^{-1}_1 (A \cup \{ a \}) \}.  
\]  
This suggests the following definition,  a measurable version of the statement of Lemma 6.13 of \cite{EGMM}.  

\begin{definition} \label{m-dis} 
Assume that a locally compact group $G$ has a strongly topological action on a topological space $\mathcal{X}$ and $g_1 \in G$. 
We say that $g_1$ is measure discerning (m-discerning) if for every $g_2 \in G$, every finite $A \subset \mathcal{X}$ and every $a\in \mathcal{X}\setminus A$ almost surely all elements $g$ of the coset $g_2 \cdot \mathsf{stab}_G (A)$ have the property that 
$\{ g(a), g_1 g (a)\} \subset \mathsf{supp}(g_1 )\setminus (A \cup \{ a \})$.  
\end{definition} 
Note that when $a\notin g^{-1}_2(A\cup \mathsf{Fix}(g_1 ))$ then existence of $g$ as in the definition also follows from the assumption of hereditary separation. 

The following proposition is a measurable version of Theorem 6.10 from \cite{EGMM}. 

\begin{proposition} \label{as}  
Assume that a locally compact group $G$ has a strongly topological 
action on a Hausdorff topological space $\mathcal{X}$. 
Then for every non-trivial word $w(y)$ which constants are m-discerning, a random element $g\in G$ almost surely satisfies the inequality $w(y) \not= 1$.  
\end{proposition} 

{\em Proof.} 
We adapt the proof of Theorem 6.10 from \cite{EGMM}. 
Assume that $w(y)$ is in the form $(4.1)$.  
Since the action is strongly topological, for a random $g\in G$ almost surely the elements $a, g(a), \ldots , g^{m_1 - 1}(a)$ are pairwise distinct (where, say, $m_1 >0$). 
Applying  m-discerning (for $g_1 = h_1$ and $g_2 = 1$) together with Fubini's theorem we see that almost surely the elements 
\[ 
a, g(a), \ldots , g^{m_1 - 1}(a), g^{m_1}(a),h_1 (g^{m_1} (a)) 
\]  
are pairwise distinct. 

Repeating the argument $k-1$ more times we obtain that almost surely the elements 
\[ 
a, g(a), \ldots , g^{m_1 - 1}(a), g^{m_1}(a),h_1 (g^{m_1} (a)), g(h_1 (g^{m_1} (a))) , \ldots , 
\] 
\[ 
\hspace{8cm} h_k ( \ldots (g^{m_2}(h_1 (g^{m_1} (a))))\ldots )
\]  
are pairwise distinct.   
We see that almost surely a random element of $G$ satisfies the inequality $w(y) \not= 1$.
$\Box$ 

\bigskip 

When $G$ acts on $\mathcal{X}$ strongly topologically and transitively,  
the space $\mathcal{X}$ can be considered under the measure inherited from $G$ (after identification of $\mathcal{X}$ with  $G/\mathsf{stab}_G (p)$). 
Then note that for each m-discerning $h\in G$ the set $\mathsf{Fix}(h)$ is of measure 0. 
In particular, in typical situations where $\mathcal{X}$ is not discrete and the conditions of Proposition \ref{as} hold,  there are non-m-discerning elements in $G$. 
For example this is the case of Corollary 1.6 of \cite{A} which concerns the automorphism group of a $d$-regular tree ($d>2$) acting on the boundary of the tree. 
On the other hand the following general statement holds. 
\begin{quote} 
{\em Assume that a locally compact group $G$ has a strongly topological 
action on a Hausdorff topological space $\mathcal{X}$.  
If every non-trivial element of $G$ is m-discerning then the group $G$ is MIF.  }
\end{quote}  
Indeed, by Remark 5.1 of \cite{HO} a group is MIF exactly when it does not have laws with constants depending on a single variable. 
It remains to apply Proposition \ref{as}. 

The discussion above suggests that the condition that all non-trivial elements of $G$ are m-discerning is very restrictive. 
The authors do not know any interesting example of this kind (in particular when $\mathcal{X}$ is not discrete and the stabilizers of finite subsets are non-trivial). 

\bigskip 

\noindent
{\bf Acknowledgement.} The authors 
are grateful to the referee for corrections and for helpful and stimulating remarks.

Aleksander Iwanow

Institute of Computer Science, University of Opole, 

ul. Oleska 48, 45 - 052 Opole, Poland 

aleksander.iwanow@uni.opole.pl 

and 
  
Department of Applied Mathematics, Silesian University of Technology, 

ul. Kaszubska 23, 44 -101 Gliwice, Poland 

Aleksander.Iwanow@polsl.pl

\bigskip 

Roland Zarzycki

Collegium Civitas, plac Defilad 1, 00-901 Warsaw, Poland 

roland.zarzycki@civitas.edu.pl

\end{document}